\documentclass[12pt,oneside,english,frenchb]{amsart}
\usepackage[latin1]{inputenc}
\usepackage{a4}
\usepackage{amsmath}
\usepackage{setspace}
\usepackage{amssymb}

\makeatletter
 \theoremstyle{plain}    
 \newtheorem{thm}{Théorème}[section]
 \numberwithin{equation}{section} 
 \numberwithin{figure}{section} 
 \theoremstyle{plain}    
 \newtheorem*{thm*}{Théorème} 
 \theoremstyle{plain}    
 \newtheorem{lem}[thm]{Lemme} 
 \theoremstyle{remark}
 \newtheorem{rem}[thm]{Remarque}
 \theoremstyle{remark}    
 \newtheorem*{claim*}{Affirmation}
 \theoremstyle{definition}
 \newtheorem{defn}[thm]{Définition}
 \theoremstyle{definition}
  \newtheorem{example}[thm]{Exemple}
 \theoremstyle{plain}    
 \newtheorem{cor}[thm]{Corollaire} 
 \theoremstyle{remark}    
 \newtheorem{notation}[thm]{Notation} 

\usepackage{calrsfs}
\usepackage{aeguill}
\usepackage{hyperref}

\DeclareMathOperator{\tr}{tr}
\DeclareMathOperator{\scal}{scal}
\DeclareMathOperator{\Ric}{Ric}

\renewcommand{\geq}{\geqslant}
\renewcommand{\leq}{\leqslant}

\usepackage{babel}
\makeatother
\begin{document}
\NoAutoSpaceBeforeFDP

\title{Autodual Einstein versus Kähler-Einstein}

\author{Olivier Biquard}

\address{Irma, Université Louis Pasteur et cnrs, 7 rue René Descartes, F-67084
Strasbourg Cedex}

\email{olivier.biquard@math.u-strasbg.fr}

\thanks{The author is a member of \textsc{edge}, Research Training Network
HPRN-CT-2000-00101, supported by the European Human Potential Programme}

\date{\today}

\maketitle

\newcommand{\bC}{\mathbb{C}}

\newcommand{\bR}{\mathbb{R}}

\newcommand{\su}{\mathfrak{su}}

\newcommand{\so}{\mathfrak{so}}

\newcommand{\slii}{\mathfrak{sl}_{2}}

\newcommand{\fm}{\mathfrak{m}}

\newcommand{\fu}{\mathfrak{u}}

\newcommand{\fC}{\mathfrak{C}}

\newcommand{\fE}{\mathfrak{E}}

\newcommand{\fS}{\mathfrak{S}}

\newcommand{\Rrond}{\overset\circ R}

\newcommand{\diagd}[2]{\begin{pmatrix}#1\\
 &  #2\end{pmatrix}}

\newcommand{\mattt}[9]{\begin{pmatrix}#1  &  #2  &  #3\\
#4  &  #5  &  #6\\
#7  &  #8  &  #9\end{pmatrix}}

\selectlanguage{english}
\begin{abstract}
Any strictly pseudoconvex domain in $\bC^{2}$ carries a complete
Kähler-Einstein metric, the Cheng-Yau metric, with ``conformal infinity''
the CR structure of the boundary.

It is well known that not all CR structures on $S^{3}$ arise in this
way. In this paper, we study CR structures on the 3-sphere satisfying
a different filling condition: boundaries at infinity of (complete)
selfdual Einstein metrics. We prove that (modulo contactomorphisms)
they form an infinite dimensional manifold, transverse to the space
of CR structures which are boundaries of complex domains (and therefore
of Kähler-Einstein metrics).
\end{abstract}
\selectlanguage{frenchb}
Un domaine strictement pseudoconvexe $D$ de $\bC^{2}$ porte une
unique métrique Kähler-Einstein complète $g$, la métrique de Cheng-Yau
\cite{CheYau80} (dans le cas de la boule $B^{4}$, il s'agit de la
métrique hyperbolique complexe, dans le modèle de la métrique de Bergmann).
Les directions complexes dans le bord $\partial D$ forment une structure
de contact, munie de la structure complexe $J$ : une structure CR
sur $\partial D$. La métrique de Cheng-Yau et la structure CR du
bord sont compatibles dans le sens suivant : choisissons une équation
$\varphi=0$ du bord, alors $\eta=Jd\varphi$ est une forme de contact
sur $\partial D$, définissant une métrique $\gamma(\cdot,\cdot)=d\eta(\cdot,J\cdot)$
sur les directions de contact, et $g$ a le comportement asymptotique\[
g\sim\frac{d\varphi^{2}+\eta^{2}}{\varphi^{2}}+\frac{\gamma}{\varphi}.\]
La métrique $\gamma$ est obtenue comme limite de $\varphi g$ sur
les directions de contact, mais n'est intrinsèque qu'à un facteur
conforme près : sa classe conforme (équivalente à la donnée de $J$)
est appelée l'\emph{infini conforme} de $g$.

Comme il est bien connu, toutes les structures CR sur $S^{3}$ ne
sont pas des bords de domaines complexes. Dans le cas où la structure
CR est proche de la structure standard, on peut cependant toujours
construire sur la boule $B^{4}$ une métrique d'Einstein, d'infini
conforme $J$ \cite{Biq00,Biq99}. Cette métrique n'est Kähler-Einstein
que dans le cas où $J$ provient d'une déformation de $S^{3}$ dans
$\bC^{2}$.

D'un autre côté, dans le cas où la structure CR $J$ est invariante
à droite, Hitchin \cite{Hit95} a démontré que $J$ est l'infini conforme
d'une métrique d'Einstein autoduale sur $B^{4}$, pour laquelle il
a fourni des formules explicites.

Dans cet article nous déterminons les structures CR, proches de la
structure standard de $S^{3}$, qui sont les infinis conformes de
métriques d'Einstein autoduales.

Bien entendu, cette condition est orthogonale à celle du remplissage
par une métrique Kähler-Einstein : en effet, si le remplissage Einstein
est à la fois Kähler et autodual, alors il est forcément hyperbolique
complexe, donc la structure CR est égal à la structure standard, à
un contactomorphisme de $S^{3}$ près.

Nous montrons que les structures CR remplissables par une métrique
autoduale Einstein remplissent toutes les directions laissées libres
par celles qui proviennent de bord de domaines complexes :

\begin{thm*}
Soit $\mathcal{C}$ l'espace des structures CR sur $S^{3}$, notons
$\mathcal{K}\subset\mathcal{C}$ celles qui sont remplissables par
une métrique Kähler-Einstein, et $\mathcal{A}\subset\mathcal{C}$
par une métrique autoduale Einstein. Enfin désignons par $J_{0}$
la structure CR standard de $S^{3}$.

Alors près de la structure standard, $\mathcal{K}$ et $\mathcal{A}$
sont deux sous-variétés transverses (leurs espaces tangents engendrent
tout l'espace tangent de $\mathcal{C}$ en $J_{0}$), et $\mathcal{K}\cap A$
est réduit à l'orbite de $J_{0}$ sous les contactomorphismes.
\end{thm*}
Le théorème est énoncé sous une forme volontairement imprécise quant
à la régularité des structures CR sur $S^{3}$ : le texte précise
les espaces fonctionnels idoines---des espaces de Folland-Stein.

L'étude de la sous-variété $\mathcal{K}$, et notamment de son espace
tangent en $J_{0}$, est dûe à Bland \cite{Bla94}. Le théorème fournit
donc seulement une assertion sur $\mathcal{A}$.

En particulier, le théorème construit une famille de dimension infinie
de nouvelles métriques autoduales Einstein ; l'espace tangent à $\mathcal{A}$
est décrit précisément dans le théorème \ref{thm:tangent}.

Le théorème est à rapprocher de la conjecture de fréquence positive
de LeBrun \cite{LeB91}, résolue dans \cite{Biq02}. Celle-ci concerne
un problème analogue, pour les métriques d'Einstein sur la boule $B^{4}$,
dont l'infini conforme est maintenant une vraie métrique conforme
sur $S^{3}$---le prototype étant fourni par la métrique hyperbolique
réelle sur la boule, induisant à l'infini la structure conforme standard
de $S^{3}$ : toute métrique proche $h$ sur $S^{3}$ se décompose
en $h=h_{+}+h_{0}+h_{-}$, où $h_{0}+h_{+}$ est l'infini conforme
d'une métrique autoduale Einstein, et $h_{0}+h_{-}$ est l'infini
conforme d'une métrique antiautoduale Einstein. 

Le théorème démontré dans ce papier indique un phénomène analogue
pour les structures CR, mais la condition d'antiautodualité (annulation
de la partie autoduale du tenseur de Weyl $W^{+}$) doit être remplacée,
car pour la métrique hyperbolique complexe (la métrique de Bergmann)
sur $B^{4}$, on a $W^{-}=0$ mais $W^{+}$, vu comme endomorphisme
des 2-formes autoduales \[
\Omega^{+}=\bR\omega\oplus\textrm{Re}(\Omega^{2,0}\oplus\Omega^{0,2}),\]
où $\omega$ est la forme de Kähler, s'écrit (avec $\textrm{scal}$
la courbure scalaire) \begin{equation}
W^{+}=\frac{\scal}{6}1_{\bR\omega}-\frac{\scal}{12}1_{\Omega^{+}}.\tag{\dag}\label{eq:Wkahler}\end{equation}
La condition de trouver une métrique à $W^{+}=0$ dans le cas réel
doit donc être modifiée dans le cas complexe : d'une certaine manière,
la condition Kähler-Einstein est la plus naturelle à substituer, car
elle consiste à figer $W^{+}$ sous la forme (\ref{eq:Wkahler}) :
un tenseur parallèle, de valeurs propres $\frac{\scal}{6}$ et $-\frac{\scal}{12}$.

Cette remarque suggère que les deux énoncés de \og fréquences positives
 \fg{}, dans les cas réel et complexe, pourraient n'être que la partie
émergée d'un phénomène plus général.

Disons quelques mots sur la démonstration : elle consiste à étudier
l'opérateur qui, à une structure CR $J$ sur le bord, associe la partie
antiautoduale $W_{g}^{-}$ du tenseur de Weyl de la métrique d'Einstein
$g$ d'infini conforme $J$. Le tenseur $W_{g}^{-}$, interprété comme
une section d'un fibré de spineurs, est harmonique car $g$ est d'Einstein.
Son comportement a été étudié dans \cite{BiqHer} : il est $L^{2}$,
et plus précisément décroît en $O(e^{-4r})$ à l'infini ($r$ étant
la distance à un point). 

L'idée principale est de montrer que l'opérateur $J\rightarrow W_{g}^{-}$
est submersif. On peut se ramener à un opérateur sur le bord, car
un tel spineur harmonique $L^{2}$ est déterminé par sa \og valeur
à l'infini \fg{}, \[
\partial W_{g}^{-}=\lim_{r\rightarrow+\infty}e^{4r}W_{g}^{-},\]
qui est une section sur $S^{3}$ du fibré $\mathcal{J}$ des endomorphismes
symétriques sans trace de la distribution de contact. Ainsi est-il
équivalent de considérer l'opérateur $J\rightarrow\partial W_{g}^{-}$.
Si on se restreint au cas où $g$ est Kähler-Einstein, alors $\partial W_{g}^{-}$
devient formellement déterminé par la structure CR $J$ : il est égal,
à une constante multiplicative près, au tenseur de Cartan $Q(J)$
de la structure CR.

Le point central de la démonstration est un calcul de \emph{tous}
les spineurs harmoniques $L^{2}$ pour la métrique de Bergmann (ils
forment un espace de dimension infinie), et en particulier de leurs
valeurs à l'infini. Cela permet de montrer qu'ils proviennent de tenseurs
de Weyl de métriques Kähler-Einstein infinitésimales ; ainsi l'opérateur
$J\rightarrow W_{g}^{-}$, à valeurs dans les spineurs harmoniques
$L^{2}$, est submersif, même restreint aux infinis conformes de métriques
Kähler-Einstein, et on en déduit le théorème.

Bien entendu, la démonstration utilise deux théorèmes difficiles :
le théorème de Bland \cite{Bla94}, identifiant les structures CR
qui sont des bords de domaines complexes (et donc de métriques Kähler-Einstein),
et le théorème de Cheng-Lee \cite{CheLee90} étudiant les propriétés
hypoelliptiques de la linéarisation de $Q$.

Une autre difficulté, plus technique, provient de la nécessité d'appliquer
le théorème des fonctions implicites dans des espaces fonctionnels
adéquats : ainsi l'opérateur $Q$ doit-il être étudié entre les espaces
de Folland-Stein sur la sphère $S^{3}$, qui sont l'équivalent des
espaces de Sobolev en géométrie de contact. En retour, cela impose
d'utiliser des infinis conformes $J$ qui ne sont pas $C^{\infty}$,
et d'obtenir, pour les métriques d'Einstein qui les remplissent, une
régularité optimale. Cette régularité pourrait provenir de l'application
d'un calcul pseudodifférentiel, comme développé dans le contexte hyperbolique
complexe par Epstein-Melrose-Mendoza \cite{EpsMelMen91}, mais, ces
auteurs se restreignant au cas scalaire, nous donnons, dans l'esprit
de \cite{Biq00}, une démonstration élémentaire aboutissant à un énoncé
de régularité de la métrique d'Einstein (lemme \ref{lem:regg}).

Donnons enfin le plan de cet article. Dans la première section, nous
étudions l'opérateur de Dirac sur l'espace hyperbolique complexe,
et déterminons explicitement, dans une décomposition harmonique, les
valeurs à l'infini des spineurs harmoniques. Dans la seconde section,
nous mettons en place l'analyse nécessaire, ce qui permet de prouver
la régularité de la métrique d'Einstein qui remplit une structure
CR donnée, et de montrer que les noyaux $L^{2}$ de l'opérateur de
Dirac sur le fibré de spineurs qui nous intéresse forment un fibré
au-dessus d'un espace adéquat de métriques. Enfin, la troisième section
est consacrée à la démonstration proprement dite du théorème.

\section{L'opérateur de Dirac sur $\bC H^{2}$\label{sec:1}}

\subsection{Décomposition harmonique}

Nous commençons par expliquer la décomposition harmonique que nous
allons faire pour étudier l'opérateur de Dirac.

\subsubsection{L'espace hyperbolique complexe}

L'espace hyperbolique complexe de dimension (réelle) 4 peut être représenté
comme \[
\bC H^{2}=SU_{1,2}/U_{1}SU_{2},\]
avec plus précisément \[
SU_{2}=\diagd{1}{*},\quad U_{1}=\mattt{z^{-2}}{}{}{}{z}{}{}{}{z}.\]
Au niveau des algèbres de Lie, on a \[
\su_{1,2}=\fu_{1}\oplus\su_{2}\oplus\fm,\]
avec $\fm$ engendré par les quatre vecteurs $X_{0},\ldots,X_{3}$
donnés explicitement par :\[
X_{0}=\mattt{}{1}{}{1}{}{}{}{}{0},\; X_{1}=\mattt{}{-i}{}{i}{}{}{}{}{0},\; X_{2}=\mattt{}{}{1}{}{0}{}{1}{}{},\; X_{3}=\mattt{}{}{-i}{}{0}{}{i}{}{}.\]
La structure complexe de $\bC H^{2}$ envoie $X_{0}$ sur $X_{1}$
et $X_{2}$ sur $X_{3}$.

La base orthonormée $(X_{i})$ de $\fm$ a été choisie de sorte que
$X_{1}$ est vecteur propre de $\textrm{ad}(X_{0})^{2}$ pour la valeur
propre $4$, et $X_{2}$ et $X_{3}$ pour la valeur propre $1$. Posons
$\lambda_{1}=2$ et $\lambda_{2}=\lambda_{3}=1$, et pour $i\geq1$,
définissons $Y_{i}\in\fu_{1}\oplus\su_{2}$ par\[
\lambda_{i}Y_{i}=[X_{0},X_{i}],\]
alors on a aussi \[
\lambda_{i}X_{i}=[X_{0},Y_{i}].\]

La signification des $(Y_{i})$ est la suivante \cite[I.1.A]{Biq00}
: le choix d'une origine $*$ dans $\bC H^{2}$ et du rayon $\exp(rX_{0})*$
permet d'identifier chaque sphère $S_{r}$ de rayon $r$ à l'espace
homogène $S^{3}=U_{1}SU_{2}/S_{1}$, par l'application $\varpi_{r}(g)=g\exp(rX_{0})*$
pour tout $g\in U_{1}SU_{2}$. Les $(Y_{i})$ représentent une base
homogène du fibré tangent à $S^{3}$ ; sur la sphère de rayon $r$,
elle est liée à la base $(X_{i})$ à l'origine par la formule \[
\exp(-rX_{0})d\varpi_{r}(Y_{i})=-\sinh(\lambda_{i}r)X_{i}.\]

Enfin, $S^{1}$ est le stabilisateur de $X_{0}$ : \[
S^{1}=\mattt{z}{}{}{}{z}{}{}{}{z^{-2}}.\]

\subsubsection{La dérivation covariante}

Soit $\mathcal{E}$ un fibré homogène sur $\bC H^{2}$, provenant
d'une représentation $\rho_{E}$ de $U_{1}SU_{2}$ dans l'espace vectoriel
$E$. Sur la sphère $S^{3}=U_{1}SU_{2}/S^{1}$, les sections $L^{2}$
de $\mathcal{E}$ se décomposent en \begin{equation}
L^{2}(S^{3},\mathcal{E})=\oplus_{\rho}V_{\rho}\otimes(V_{\rho}\otimes E)^{S^{1}},\label{eq:L2}\end{equation}
où la somme court sur toutes les représentations irréductibles $\rho$
de $U_{1}SU_{2}$ (dans $V_{\rho}$). La section de $\mathcal{E}$
correspondant à $v\otimes w\in V_{\rho}\otimes(V_{\rho}\otimes E)^{S^{1}}$est
donnée par l'application $S^{1}$-équivariante $s_{v\otimes w}:U_{1}SU_{2}\rightarrow V_{0}$
définie par \[
s_{v\otimes w}(g)=\left\langle w,\rho(g^{-1})v\right\rangle .\]
Une section de $\mathcal{E}$ sur $\mathbb{C}H^{2}$ tout entier se
décompose ainsi sur la somme (\ref{eq:L2}) avec des coefficients
dépendant du rayon $r$ : il sera commode de prendre $v$ fixe et
$w$ dépendant de $r$, soit $s=v\otimes w(r)$.

Ces notations posées, la dérivation covariante sur $\mathcal{E}$
est donnée par les formules (voir \cite[I.2.A]{Biq00})\begin{align}
\nabla_{X_{0}}s & =v\otimes\partial_{r}w,\label{eq:nabla0}\\
\nabla_{X_{i}}s & =-v\otimes\left(\coth(\lambda_{i}r)\rho_{E}(Y_{i})+\frac{1}{\sinh(\lambda_{i}r)}\rho(Y_{i})\right)w,\quad i\geq1.\label{eq:nabla1}\end{align}
Bien entendu, les représentations $\rho_{E}$ et $\rho$ qui apparaissent
dans ces formules sont en réalité les différentielles $(\rho_{E})_{*}$
et $\rho_{*}$ induites sur l'algèbre de Lie $\fu_{1}\oplus\su_{2}$
; en l'absence d'ambiguïté, nous les notons par le même symbole.

\subsubsection{L'action des $Y_{i}$}

Donnons tout d'abord une base de $\so(\fm)=\Omega^{2}\fm$ ; il est
bien connu qu'il y a une décomposition $\so_{4}=\su_{2}\oplus\su_{2}$
correspondant à la décomposition en formes autoduales et antiautoduales
$\Omega^{2}=\Omega^{+}\oplus\Omega^{-}$ ; explicitement, si $*$
est l'opérateur de Hodge (donc $*(X_{0}\wedge X_{1})=X_{2}\wedge X_{3}$,
etc.), on peut poser pour $i\geq1$\begin{align}
\sigma_{i}^{+} & =X_{0}\wedge X_{i}+*(X_{0}\wedge X_{i})\label{eq:siplus}\\
\sigma_{i}^{-} & =-X_{0}\wedge X_{i}+*(X_{0}\wedge X_{i})\label{eq:simoins}\end{align}
Les signes ont été choisis de sorte que l'on ait les mêmes relations
pour les deux bases $(\sigma_{i}^{+})$ et $(\sigma_{i}^{-})$, à
savoir \begin{equation}
[\sigma_{1},\sigma_{2}]=2\sigma_{3},\quad[\sigma_{2},\sigma_{3}]=2\sigma_{1},\quad[\sigma_{3},\sigma_{1}]=2\sigma_{2}.\label{eq:relsu2}\end{equation}

Maintenant, on peut calculer explicitement \[
Y_{1}=\mattt{i}{}{}{}{-i}{}{}{}{0},\quad Y_{2}=\mattt{0}{}{}{}{}{1}{}{-1}{},\quad Y_{3}=\mattt{0}{}{}{}{}{-i}{}{-i}{},\]
et on vérifie aisément que l'action adjointe des $Y_{i}$ sur $\fm$
est donnée par \begin{equation}
Y_{1}=\frac{1}{2}\sigma_{1}^{-}-\frac{3}{2}\sigma_{1}^{+},\quad Y_{2}=\sigma_{2}^{-},\quad Y_{3}=\sigma_{3}^{-}.\label{eq:Yi}\end{equation}
Cela explicite la formule (\ref{eq:nabla1}).

Finalement, l'action infinitésimale du $S^{1}$ qui stabilise $X_{0}$
est donnée par \begin{equation}
\xi=\sigma_{1}^{+}+\sigma_{1}^{-}.\label{eq:xi}\end{equation}

\subsection{L'opérateur de Dirac}

\subsubsection{Les spineurs de $\bC H^{2}$ et l'opérateur de Dirac}

Soit $\mathcal{S}=\mathcal{S}^{+}\oplus\mathcal{S}^{-}$les spineurs
de $\bC H^{2}$. Comme représentations de $U_{1}SU_{2}$, le fibré
$\mathcal{S}^{+}$ provient de la représentation $\rho^{+}$ de $U_{1}$
sur $S^{+}=\bC_{1}\oplus\bC_{-1}$, tandis que $\mathcal{S}^{-}$
provient de la représentation tautologique $\rho^{-}$ de $SU_{2}$
sur $S^{-}=\bC^{2}$.

Soit $\mathcal{E}$ un fibré homogène, associé à la représentation
$\rho_{E}$ de $U_{1}SU_{2}$, alors l'opérateur de Dirac tordu $\mathcal{D}:\mathcal{S}^{-}\otimes\mathcal{E}\rightarrow\mathcal{S}^{+}\otimes\mathcal{E}$
est donné par \[
\mathcal{D}s=\sum_{1}^{4}X_{i}\cdot\nabla_{X_{i}}s.\]
En termes de la décomposition harmonique (\ref{eq:L2}), compte tenu
des formules (\ref{eq:nabla0}) et (\ref{eq:nabla1}), cela se traduit
sur la section $s=v\otimes w(r)$, où $w(r)\in(V_{\rho}\otimes S^{-}\otimes E)^{S^{1}}$,
par un opérateur agissant sur $w$ seulement : \[
\mathcal{D}w=X_{0}\cdot\partial_{r}w-\sum_{1}^{3}X_{i}\cdot\left(\coth(\lambda_{i}r)\rho_{0}(Y_{i})+\frac{1}{\sinh(\lambda_{i}r)}\rho(Y_{i})\right)w,\]
où $\rho_{0}=\rho^{-}\otimes\rho_{E}$ est la représentation de $U_{1}SU_{2}$
sur $S^{-}\otimes E$. Le résultat est une fonction de $r$, à valeurs
dans $(V_{\rho}\otimes S^{+}\otimes E)^{S^{1}}$ ; pour revenir dans
le même espace $(V_{\rho}\otimes S^{-}\otimes E)^{S^{1}}$, il sera
commode de considérer l'opérateur \[
\partial_{r}\cdot\mathcal{D}w=-\partial_{r}w-\sum_{1}^{3}X_{0}\cdot X_{i}\cdot\left(\coth(\lambda_{i}r)\rho_{0}(Y_{i})+\frac{1}{\sinh(\lambda_{i}r)}\rho(Y_{i})\right)w.\]
L'action de Clifford de $X_{0}X_{i}$ sur $S^{-}$ n'est autre que
$-\rho^{-}(\sigma_{i}^{-})$, d'où résulte la formule : \begin{equation}
\partial_{r}\cdot\mathcal{D}w=-\partial_{r}w+\sum_{1}^{3}\rho^{-}(\sigma_{i}^{-})\left(\coth(\lambda_{i}r)\rho_{0}(Y_{i})+\frac{1}{\sinh(\lambda_{i}r)}\rho(Y_{i})\right)w.\label{eq:drD}\end{equation}

\subsubsection{L'opérateur de Dirac sur $\mathcal{S}_{4}^{-}$}

On notera \[
S_{\ell}^{\pm}=\fS^{\ell}S^{\pm},\quad\mathcal{S}_{\ell}^{\pm}=\fS^{\ell}\mathcal{S}^{\pm}\]
les puissances symétriques $\ell$-ièmes. Comme $S^{-}S_{3}^{-}=S_{4}^{-}\oplus S_{2}^{-}$,
on dispose d'un opérateur de Dirac \[
\mathcal{D}:\mathcal{S}_{4}^{-}\longrightarrow\mathcal{S}^{+}\mathcal{S}_{3}^{-}.\]

\begin{lem}
\label{lem:DS4}Si la métrique est autoduale, alors $\mathcal{D}^{2}:\mathcal{S}^{-}\mathcal{S}_{3}^{-}\longrightarrow\mathcal{S}^{-}\mathcal{S}_{3}^{-}$
préserve $\mathcal{S}_{4}^{-}$ et $\mathcal{S}_{2}^{-}$.
\end{lem}
\begin{proof}
La projection $\mathcal{D}^{2}:\mathcal{S}_{4}^{-}\rightarrow\mathcal{S}_{2}^{-}$
est un opérateur tensoriel donnée par l'action d'un morceau de la
courbure ; le seul possible ici est $W^{-}\in\mathcal{S}_{4}^{-}$.
\end{proof}
L'opérateur \[
\partial_{r}\cdot\mathcal{D}:\mathcal{S}^{-}\mathcal{S}_{3}^{-}\longrightarrow\mathcal{S}^{-}\mathcal{S}_{3}^{-}\]
a été calculé par la formule (\ref{eq:drD}) : dans la somme $\mathcal{S}^{-}\mathcal{S}_{3}^{-}=\mathcal{S}_{4}^{-}\oplus\mathcal{S}_{2}^{-}$,
il se décompose en \[
\partial_{r}\cdot\mathcal{D}=\left(\begin{array}{cc}
-\partial_{r}+A & C\\
B & -\partial_{r}+D\end{array}\right),\]
où $A$, $B$, $C$ et $D$ sont des opérateurs différentiels qui
ne dérivent les spineurs que dans la direction des sphères $S_{r}$.

\begin{lem}
On a l'identité \[
\mathcal{D}\partial_{r}\cdot+\partial_{r}\cdot\mathcal{D}=-2\nabla_{\partial_{r}}+\tr\mathbb{I},\]
où $\mathbb{I}$ est la seconde forme fondamentale de la sphère $S_{r}$.
\end{lem}
\begin{proof}
On a les égalités, en utilisant notre base orthonormale $(X_{i})$
telle que $X_{0}=\partial_{r}$, \begin{align*}
\mathcal{D}\partial_{r}\cdot\sigma & =\sum_{0}^{3}X_{i}\cdot\nabla_{X_{i}}(\partial_{r}\sigma)\\
 & =\sum_{0}^{3}X_{i}\cdot(\nabla_{X_{i}}\partial_{r})\cdot\sigma+X_{i}\cdot\partial_{r}\cdot\nabla_{X_{i}}\sigma\\
 & =(\tr\mathbb{I})\sigma-\nabla_{\partial_{r}}\sigma-\sum_{1}^{3}\partial_{r}\cdot X_{i}\cdot\nabla_{X_{i}}\sigma.\end{align*}

\end{proof}
Il résulte du lemme que l'on a l'identité \[
(\partial_{r}\cdot\mathcal{D})^{2}=\mathcal{D}^{2}+(\tr\mathbb{I}-2\nabla_{\partial_{r}})\partial_{r}\cdot\mathcal{D}.\]
 Compte tenu du lemme \ref{lem:DS4}, cela impose l'identité\begin{equation}
\nabla_{\partial_{r}}B+BA=(\tr\mathbb{I}-D)B.\label{eq:drB}\end{equation}
On en déduit le résultat suivant :

\begin{lem}
\label{lem:intD}Supposons que l'on ait une section $\sigma$ de $\mathcal{S}_{4}^{-}$
telle que :
\begin{enumerate}
\item $(-\partial_{r}+A)\sigma=0$;
\item $B\sigma=0$ pour une valeur de $r$ ;
\end{enumerate}
alors $B\sigma=0$ pour tout $r$, c'est-à-dire $\mathcal{D}\sigma=0$.
\end{lem}
\begin{proof}
L'équation (\ref{eq:drB}) et la première condition impliquent \[
\nabla_{\partial_{r}}(B\sigma)=(\tr\mathbb{I}-D)(B\sigma),\]
ce qui donne le résultat.
\end{proof}

\subsubsection{Décomposition harmonique}

Analysons maintenant cet opérateur de Dirac du point de vue de la
décomposition harmonique. La décomposition (\ref{eq:L2}) fait intervenir
des éléments $w\in(V_{\rho}S_{4}^{-})^{S^{1}}$. Or, les représentations
$\rho$ de $U_{1}SU_{2}$ sont indexées par deux entiers $K\in\mathbb{Z}$
et $L\in\mathbb{N}$, de sorte que \[
V_{\rho}=\bC_{K}\otimes S^{L},\]
où $U_{1}$ agit avec poids $K$ sur $\bC_{K}$ et $SU_{2}$ agit
sur $S^{L}=S_{-}^{L}$. D'après (\ref{eq:xi}), l'action du générateur
$\xi=\sigma_{1}^{+}+\sigma_{1}^{-}$ de $S^{1}$ sur $V_{\rho}S_{4}^{-}$
est par \[
i(K+k_{1}+k_{2}),\]
où les $k_{1}\in\{-L,-L+2,\cdots,L\}$ sont les valeurs propres de
l'action de $-i\sigma_{1}^{-}$ sur $S^{L}$, et $k_{2}\in\{0,\pm2,\pm4\}$
celles de l'action sur $S_{-}^{4}$. En particulier : \begin{equation}
k_{2}=-K-k_{1},\label{eq:kkK}\end{equation}
donc $(V_{\rho}S_{4}^{-})^{S^{1}}=0$ si $|K|>L+4$, et en général
\[
\dim(V_{\rho}S_{4}^{-})^{S^{1}}\leq5,\]
avec égalité si $|K|\leq L-4$. 

On va en déduire :

\begin{lem}
\label{lem:solD0}L'espace des solutions de $\mathcal{D}\sigma=0$
sur $\bC H^{2}-\{*\}$ avec $\sigma$ provenant de la représentation
$V_{\rho}=\bC_{K}S^{L}$ comme dans (\ref{eq:L2}) est de dimension
:
\begin{itemize}
\item $2\dim V_{\rho}$ si $|K|\leq L-4$ ;
\item $\dim V_{\rho}$ si $L-4<|K|\leq L+4$.
\end{itemize}
\end{lem}
\begin{proof}
Pour le spineur $\sigma$ correspondant à $v\otimes w(r)$, où $v$
est fixe, l'équation $\mathcal{D}\sigma=0$ se réduit à l'équation
différentielle ordinaire $(-\partial_{r}+A)w=0$ sous la contrainte
algébrique $Bw=0$ :
\begin{itemize}
\item pour $|K|\leq L-4$, $w(r)\in(V_{\rho}S_{4}^{-})^{S^{1}}$ qui est
de dimension 5, donc l'espace des solutions de l'équation $(-\partial_{r}+A)w=0$
est de dimension 5 ; on doit ajouter la contrainte $Bw(r)=0$, mais
$Bw(r)\in(V_{\rho}S_{2}^{-})^{S^{1}}$ avec $\dim(V_{\rho}S_{2}^{-})^{S^{1}}=3$
: compte tenu du lemme \ref{lem:intD}, il suffit d'imposer la contrainte
pour un seul $r$ et donc, $v$ étant fixé, l'espace des solutions
$w$ telles que $\mathcal{D}\sigma=0$ est de dimension 2 ;
\item pour $L-4<|K|\leq L+4$, la même analyse de dimension amène à un espace
de solutions de dimension 1.
\end{itemize}
\end{proof}

\subsection{Spineurs harmoniques $L^{2}$}

Le lemme \ref{lem:solD0} donne l'espace des solutions de l'équation
$\mathcal{D}\sigma=0$ en dehors de l'origine $*$. Nous regardons
maintenant quelles solutions se prolongent à l'origine. Puisque $\mathcal{D}$
est un opérateur elliptique, cela est le cas si et seulement si $\sigma$
est dans l'espace $L^{2}$ au voisinage de $*$. Nous déterminons
également les solutions $L^{2}$ à l'infini.

\subsubsection{Calcul de l'opérateur de Dirac sur $\mathcal{S}_{4}^{-}$}

Nous appliquons la formule (\ref{eq:drD}) pour l'opérateur de Dirac
sur $\mathcal{S}_{4}^{-}$, donc $\rho_{0}$ est maintenant la représentation
sur $S_{4}^{-}$, en particulier $\rho_{0}(\sigma_{1}^{+})=0$. Compte
tenu de $Y_{1}=\frac{1}{2}\sigma_{1}^{-}-\frac{3}{2}\sigma_{1}^{+}$
par (\ref{eq:Yi}), on obtient\begin{align*}
\partial_{r}\cdot\mathcal{D} & =-\partial_{r}+\coth r\sum_{1}^{3}\rho^{-}(\sigma_{i}^{-})\rho_{0}(\sigma_{i}^{-})+\frac{1}{\sinh r}\sum_{1}^{3}\rho^{-}(\sigma_{i}^{-})\rho(\sigma_{i}^{-})\\
 & \qquad+\rho^{-}(\sigma_{1}^{-})\left((\frac{\coth2r}{2}-\coth r)\rho_{0}(\sigma_{1}^{-})+(\frac{1}{2\sinh2r}-\frac{1}{\sinh r})\rho(\sigma_{1}^{-})\right)\\
 & \qquad-\frac{3}{2\sinh2r}\rho^{-}(\sigma_{1}^{-})\rho(\sigma_{1}^{+})\end{align*}
Puisque le générateur du stabilisateur $S^{1}$ est $\sigma_{1}^{-}+\sigma_{1}^{+}$,
l'action sur $(V_{\rho}\otimes S_{4}^{-})^{S^{1}}$ de $(\rho+\rho_{0})(\sigma_{1}^{+})=\rho(\sigma_{1}^{+})$
est égale à celle de $-(\rho+\rho_{0})(\sigma_{1}^{-})$, et on obtient
par un calcul facile \begin{align}
\partial_{r}\cdot\mathcal{D} & =-\partial_{r}+\coth r\sum_{1}^{3}\rho^{-}(\sigma_{i}^{-})\rho_{0}(\sigma_{i}^{-})+\frac{1}{\sinh r}\sum_{1}^{3}\rho^{-}(\sigma_{i}^{-})\rho(\sigma_{i}^{-})\label{eq:drDbis}\\
 & \qquad-\frac{\tanh r}{2}\rho^{-}(\sigma_{1}^{-})\rho_{0}(\sigma_{1}^{-})-\frac{\tanh r}{1+\cosh r}\rho^{-}(\sigma_{1}^{-})\rho(\sigma_{1}^{-})\,.\nonumber \end{align}

\subsubsection{Comportement à l'origine}

Poursuivons le calcul quand $r\rightarrow0$, on obtient \[
\partial_{r}\cdot\mathcal{D}w\sim-\partial_{r}w+\frac{1}{r}\sum_{1}^{3}\rho^{-}(\sigma_{i}^{-})(\rho_{0}+\rho)(\sigma_{i}^{-})w.\]
 Voyant $\rho_{0}$ comme la représentation sur $S_{4}^{-}\subset S^{-}S_{3}^{-}$,
et en tenant compte de l'identité $\sum_{1}^{3}\rho^{-}(\sigma_{i}^{-})^{2}=-3$,
il suit\begin{equation}
\partial_{r}\cdot\mathcal{D}w\sim-\partial_{r}w-\frac{3w}{r}+\frac{1}{r}\sum_{1}^{3}\rho^{-}(\sigma_{i}^{-})(\rho_{3}^{-}+\rho)(\sigma_{i}^{-})w.\label{eq:drD0}\end{equation}

\begin{lem}
\label{lem:cas}Sur la représentation $S\otimes S_{\ell}=S_{\ell+1}\oplus S_{\ell-1}$
de $\su_{2}$, on a l'identité \[
\sum_{1}^{3}\rho_{1}(\sigma_{i})\rho_{\ell}(\sigma_{i})=\begin{cases}
-\ell & \textrm{ sur }S_{\ell+1}\\
\ell+2 & \textrm{ sur }S_{\ell-1}\end{cases}\]

\end{lem}
\begin{proof}
Tout d'abord, l'opérateur de Casimir sur $S_{\ell}$ est donné par
\[
\fC(\rho_{\ell})=-\sum_{1}^{3}\rho_{\ell}(\sigma_{i})^{2}=\ell(\ell+2).\]
Puis on écrit \begin{align*}
\sum_{1}^{3}\rho_{1}(\sigma_{i})\rho_{\ell}(\sigma_{i}) & =\frac{1}{2}\left(\sum_{1}^{3}(\rho_{1}+\rho_{\ell})(\sigma_{i})^{2}-\rho_{1}(\sigma_{i})^{2}-\rho_{\ell}(\sigma_{i})^{2}\right)\\
 & =\frac{1}{2}\left(-\fC(\rho_{1}\otimes\rho_{\ell})+\fC(\rho_{1})+\fC(\rho_{\ell})\right)\end{align*}
d'où on déduit le résultat, vu la décomposition $\rho_{1}\otimes\rho_{\ell}=\rho_{\ell+1}\oplus\rho_{\ell-1}$.
\end{proof}
Calculons à présent le terme algébrique dans l'équation (\ref{eq:drD0})
: nous avons une décomposition \[
S_{3}S_{L}=S_{L+3}\oplus S_{L+1}\oplus S_{L-1}\oplus S_{L-3}.\]
Tensorisant avec la représentation $S=S_{1}$, les composantes de
$S_{1}S_{3}S_{L}$ correspondantes sont données par le tableau de
représentations :\begin{equation}
\begin{array}{cccc}
S_{L+4} & S_{L+2} & S_{L} & S_{L-2}\\
S_{L+2} & S_{L} & S_{L-2} & S_{L-4}\end{array}\label{eq:decSSS}\end{equation}
En appliquant le lemme \ref{lem:cas}, on obtient le tableau suivant
de valeurs propres pour le terme $\sum_{1}^{3}\rho_{1}(\sigma_{i})(\rho_{3}+\rho_{L})(\sigma_{i})$
: \begin{equation}
\begin{array}{cccc}
-(L+3) & -(L+1) & -(L-1) & -(L-3)\\
L+5 & L+3 & L+1 & L-1\end{array}\label{eq:vpSSS}\end{equation}

Près de l'origine, l'équation (\ref{eq:drD0}) nous indique que, sur
chaque morceau de la décomposition de $S_{1}S_{3}S_{L}$ dans (\ref{eq:decSSS}),\begin{equation}
\partial_{r}\cdot\mathcal{D}w\sim(-\partial_{r}+\frac{-3+\lambda}{r})w,\label{eq:cmpt0}\end{equation}
où les $\lambda$ sont les valeurs propres calculées en (\ref{eq:vpSSS}).
Cela nous fournit des solutions de Dirac avec comportement asymptotique
en $r^{-3+\lambda}$ près de l'origine, donc les solutions qui se
prolongent à l'origine, c'est-à-dire qui sont $L^{2}$ près de l'origine,
sont celles avec comportement asymptotique dans les composantes à
\[
\lambda>\frac{3}{2}.\]

De plus, nous voulons nous limiter aux solutions $s\in\mathcal{S}_{4}^{-}$.
Regardons d'abord le cas où $|K|\leq L-4$, donc $w(r)\in(V_{\rho}\otimes S_{4}^{-})^{S^{1}}$
qui est de dimension $5$. Remarquons que, dans le tableau (\ref{eq:decSSS}),
les représentations $S_{L\pm4}$ ne peuvent provenir que de $S_{L}S_{4}\subset S_{L}S_{3}S_{1}$
(et pas de $S_{L}S_{2}\subset S_{L}S_{3}S_{1}$), alors que les autres
représentations peuvent provenir des deux morceaux. Algébriquement,
cela signifie que $(V_{\rho}\otimes S_{4}^{-})^{S^{1}}$ et $(V_{\rho}\otimes S_{4}^{-})^{S^{1}}$
se décomposent sur (\ref{eq:decSSS}) en \begin{align*}
(V_{\rho}\otimes S_{4}^{-})^{S^{1}} & =P_{L+4}\oplus P_{L+2}\oplus P_{L}\oplus P_{L-2}\oplus P_{L-4},\\
(V_{\rho}\otimes S_{2}^{-})^{S^{1}} & =Q_{L+2}\oplus Q_{L}\oplus Q_{L-2},\end{align*}
avec tous les $P_{i}$ et les $Q_{i}$ de dimension 1, et $P_{L\pm4}\subset S_{L\pm4}$,
mais pour $i=0,\pm2$, \[
P_{L\pm i}\oplus Q_{L\pm i}\subset S_{L\pm i}\oplus S_{L\pm i}.\]
Quand $r\rightarrow0$, la forme asymptotique (\ref{eq:cmpt0}) de
l'opérateur de Dirac $\partial_{r}\cdot\mathcal{D}w(r)$ préserve
donc $P_{L\pm4}$. En revanche, pour $i=0,\pm2$, les valeurs propres
$\lambda$ sur les deux copies de $S_{L\pm i}$, calculées en (\ref{eq:vpSSS}),
sont distinctes, donc $P_{L\pm i}$ n'est pas préservé : la projection
sur $Q_{L\pm i}$ est non triviale, c'est-à-dire la projection sur
$\mathcal{S}_{2}^{-}$ de $\partial_{r}\cdot\mathcal{D}s$ est non
nulle. On en déduit que l'espace des $w(r)\in(S_{4}^{-}V_{\rho})^{S^{1}}$
tels que $\mathcal{D}w=0$, qui est de dimension $2$ d'après le lemme
\ref{lem:solD0}, possède quand $r\rightarrow0$ un comportement asymptotique
donné par des $w(r)\in P_{L+4}\oplus P_{L-4}$ (et réciproquement,
un tel comportement asymptotique détermine une unique solution). Le
calcul des valeurs propres $\lambda$ nous indique que la condition
$\lambda>-1$ n'est vérifiée que sur $P_{L-4}$. Ainsi l'espace des
solutions de $\mathcal{D}\sigma=0$ provenant de la représentation
$\rho$, et se prolongeant à l'origine, est de dimension $\dim V_{\rho}$.
On en déduit la première partie du lemme suivant :

\begin{lem}
\label{lem:solD1}L'espace des solutions de $\mathcal{D}s=0$ sur
$\bC H^{2}$ avec $s$ provenant de la représentation $V_{\rho}=\bC_{K}S^{L}$
comme dans (\ref{eq:L2}) est de dimension :
\begin{itemize}
\item $\dim V_{\rho}$ si $|K|\leq L-4$ ;
\item $0$ si $L-4<|K|\leq L+4$.
\end{itemize}
\end{lem}
\begin{proof}
Il ne reste à analyser que le cas $L-4<|K|\leq L+4$ : de manière
similaire, seule la composante $P_{L-4}$ peut fournir des solutions
$L^{2}$ pour l'opérateur de Dirac, donc il suffit de voir quand elle
est encore présente. D'après l'égalité (\ref{eq:kkK}), les poids
$k_{1}$ sur $S_{L}$ et $k_{2}$ sur $S_{4}$ subissent la contrainte
$k_{1}+k_{2}=-K$ ; or le poids correspondant sur $S_{L}S_{4}$ est
justement $k_{1}+k_{2}$, donc une projection non nulle sur $S_{L-4}$
ne peut exister que si $|k_{1}+k_{2}|\leq L-4$, soit $|K|\leq L-4$.
\end{proof}

\subsubsection{Comportement à l'infini}

Considérons à nouveau $s$ provenant de $v\otimes w(r)\in V_{\rho}\otimes(V_{\rho}\otimes S_{4}^{-})^{S^{1}}$
pour la représentation $\rho$ de $U_{1}SU_{2}$ sur $V_{\rho}=\bC_{K}\otimes S^{L}$.
À partir de l'équation (\ref{eq:drDbis}), avec $\rho_{0}$ la représentation
sur $S_{4}^{-}$, on obtient : \[
\partial_{r}\cdot\mathcal{D}w\sim-\partial_{r}w-\frac{1}{2}\rho^{-}(\sigma_{1}^{-})\rho_{0}(\sigma_{1}^{-})w+\sum_{1}^{3}\rho^{-}(\sigma_{i}^{-})\rho_{0}(\sigma_{i}^{-})w.\]
Voyant $S_{4}^{-}\subset S^{-}S_{3}^{-}$, cela nous donne \[
\partial_{r}\cdot\mathcal{D}w\sim-\partial_{r}w-\frac{5}{2}w-\frac{1}{2}\rho^{-}(\sigma_{1}^{-})\rho_{3}^{-}(\sigma_{1}^{-})w+\sum_{1}^{3}\rho^{-}(\sigma_{i}^{-})\rho_{3}^{-}(\sigma_{i}^{-})w\]
Le dernier terme a été calculé dans le lemme \ref{lem:cas}, et est
égal à $-3$ sur $S_{4}^{-}$. On en déduit \[
\partial_{r}\cdot\mathcal{D}w\sim-\partial_{r}w-\frac{11}{2}w-\frac{1}{2}\rho^{-}(\sigma_{1}^{-})\rho_{3}^{-}(\sigma_{1}^{-})w.\]

Les valeurs propres de $\rho^{-}(\sigma_{1}^{-})$ sont $\pm i$,
et celles de $\rho_{3}^{-}(\sigma_{1}^{-})$ sont $\pm i$ et $\pm3i$.
On en déduit que la plus grande valeur propre de $-\rho^{-}(\sigma_{1}^{-})\rho_{3}^{-}(\sigma_{1}^{-})$
est $3$, et par conséquent on obtient $w=O(e^{-4r})$, donc les solutions
sont $L^{2}$ à l'infini ; plus précisément, on peut même écrire \begin{equation}
w=w_{\infty}e^{-4r}+O(e^{-5r}),\label{eq:winf}\end{equation}
où $w_{\infty}$ est dans l'espace propre de $-\rho^{-}(\sigma_{1}^{-})\rho_{3}^{-}(\sigma_{1}^{-})$
pour la valeur propre $3$. Cet espace propre n'est autre que le sous-espace
$P_{-4}\oplus P_{4}\subset S_{4}^{-}$ des vecteurs de poids $\pm4$,
il est donc de dimension $2$.

On a ainsi obtenu le résultat suivant.

\begin{lem}
\label{lem:ds}Toutes les solutions $s$ de $\mathcal{D}s=0$, provenant
de la représentation $\rho$, sont $L^{2}$ à l'infini. Plus précisément,
elles vérifient \[
s=s_{\infty}e^{-4r}+O(e^{-5r}),\]
où $s_{\infty}$ est une section sur la sphère à l'infini $S^{3}$
du sous-fibré $\mathcal{J}$ de $\mathcal{S}_{4}^{-}$ constitué des
spineurs de poids $\pm4i$ pour l'action de $\sigma_{1}^{-}$.
\end{lem}
\begin{rem}
La définition (\ref{eq:simoins}) de $\sigma_{1}^{-}$ est intrinsèque
à l'infini, puisque $X_{0}=\partial_{r}$ et $X_{1}=JX_{0}$. En identifiant
$\mathcal{S}_{-}^{4}$ avec $\fS_{0}^{2}\Omega^{-}$, on peut identifier
$\mathcal{J}$ aux formes engendrées par $\sigma_{2}^{-}\sigma_{2}^{-}-\sigma_{3}^{-}\sigma_{3}^{-}$
et $\sigma_{2}^{-}\sigma_{3}^{-}+\sigma_{3}^{-}\sigma_{2}^{-}$.

Autre interprétation : l'application $\omega\rightarrow-X_{0}\lrcorner\omega$
permet d'identifier $\Omega^{-}$ avec l'espace engendré par $(X_{1},X_{2},X_{3})$,
et le sous-fibré $\mathcal{J}$ est l'espace des endomorphismes symétriques,
sans trace, de l'espace engendré par $X_{2}$ et $X_{3}$. Sur le
bord, on a $TS^{3}=\bR R\oplus\ker\eta$, avec $R$ le champ de Reeb
de la forme de contact standard $\eta$ sur $S^{3}$, et la base $(X_{0},\dots,X_{3})$
s'identifie à \[
(\partial_{r},\sinh(2r)R,\sinh(r)h,\sinh(r)J_{0}h),\]
 où $h$ est un vecteur unitaire de $\ker\eta$. Ainsi, le sous-fibré
$\mathcal{J}$ s'identifie à $\fS_{0}^{2}\ker\eta$, c'est-à-dire
au fibré des déformations infinitésimales de la structure CR, voir
(\ref{eq:phi}). Il est important de noter le poids qui apparaît :
une section $\Gamma$ de $\fS_{0}^{2}\ker\eta$ donne la section $\sinh^{2}(r)\Gamma$
de $\mathcal{J}$, qui est de norme non nulle par rapport à la métrique
hyperbolique complexe.
\end{rem}

\begin{rem}
\label{rem:obstruction}Les sections $s_{\infty}$ de $\mathcal{J}$
ne sont pas toutes des valeurs à l'infini de spineurs harmoniques,
même localement : en effet, on peut calculer que la résolution formelle
du problème $\mathcal{D}s=0$ avec $s_{\infty}$ donnée admet une
obstruction donnée par un opérateur différentiel d'ordre 2 sur $s_{\infty}$. 
\end{rem}

\begin{rem}
L'application \og valeur à l'infini \fg{} fournie par le lemme est
donnée entre des espaces de dimension finie, provenant de la représentation
$\rho$. Le lemme ne dit rien sur la régularité de la valeur au bord
d'un spineur harmonique quelconque.
\end{rem}

\subsection{Calcul et conclusion}

Les solutions de l'équation $\mathcal{D}s=0$ provenant de la représentation
$\rho$ forment un espace de dimension $2$ d'après le lemme \ref{lem:solD0}.
Par conséquent, toujours en se restreignant à $\rho$, l'application
valeur à l'infini, $s\rightarrow s_{\infty}$, fournie par le lemme
\ref{lem:ds}, est une application entre deux espaces de dimension
$2$. Par un calcul explicite, on va maintenant montrer qu'il s'agit
en fait d'un isomorphisme, et donc que l'application qui à un spineur
harmonique $L^{2}$, $s$, associe sa valeur à l'infini $s_{\infty}$,
est injective.

\subsubsection{Solutions explicites de l'équation de Dirac}

Nous ramenons l'équation $\mathcal{D}s=0$, qui est, sur chaque représentation,
un système différentiel d'ordre 1, à une seule équation d'ordre 2,
dont les solutions sont fournies en terme de fonctions hypergéométriques.

Pour alléger la rédaction, certaines étapes intermédiaires seront
laissées au lecteur.

Rappelons qu'à partir de la représentation de $\su_{2}$ donnée par
les relations (\ref{eq:relsu2}), on obtient une représentation $(H,X,Y)$
de $\slii$ en posant $\sigma_{1}=iH$, $Y=\frac{1}{2}(\sigma_{2}+i\sigma_{3})$,
$X=\frac{1}{2}(-\sigma_{2}+i\sigma_{3})$. Ainsi, $[H,Y]=-2Y$, $[H,X]=2X$
et $[X,Y]=H$.

Revenons à l'opérateur de Dirac (\ref{eq:drDbis}), et écrivons, pour
une section de $\mathcal{S}_{4}^{-}$, \[
\partial_{r}\cdot\mathcal{D}=-\partial_{r}-6\coth r+\frac{\tanh r}{2}+\frac{\tanh r}{2}\, A-\frac{2}{\sinh2r}\, B-\frac{2}{\sinh r}\, C,\]
où $A$, $B$, et $C$ sont les opérateurs algébriques donnés par
\[
A=\rho^{-}(H)\rho_{3}^{-}(H),\qquad B=\rho^{-}(H)\rho(H),\]
\[
C=-\frac{1}{2}\sum_{2}^{3}\rho^{-}(\sigma_{i}^{-})\rho(\sigma_{i}^{-})=\rho^{-}(X)\rho(Y)+\rho^{-}(Y)\rho(X).\]
Il s'avérera utile de faire le changement de variable $u=\sinh^{2}r$,
de sorte que l'opérateur à analyser devient \[
P=-2u(1+u)\partial_{u}-6-\frac{11}{2}u+\frac{u}{2}A-B-2\sqrt{1+u}C.\]

Pour calculer ces opérateurs, nous choisissons des bases construites
de la manière suivante : soit $e$ un vecteur de plus haut poids pour
$S$ (ainsi $(e,Ye)$ forme une base de $S$), $f$ un vecteur de
plus haut poids pour $S_{3}$ (donc $(e,\dots,Y^{3}e)$ est une base
de $S_{3}$) et $g$ un vecteur de poids $K$ dans $S_{L}$.

Des vecteurs de plus haut poids $s_{4}$ et $t_{2}$ pour $S_{4}$
et $S_{2}$ vus à l'intérieur de $SS_{3}$ sont donnés par \[
s_{4}=ef,\quad t_{2}=e(Yf)-3(Ye)f.\]
À partir de là, on obtient pour $S_{4}$ et $S_{2}$ les bases $(s_{4},s_{2}=Ys_{4},s_{0}=Y^{2}s_{4},\dots)$
et $(t_{2},t_{0}=Yt_{2},t_{-2}=Y^{2}t_{2})$.

Plaçons-nous dans le cas où $|K|\leq L-4$, ce qui est légitime puisque,
d'après le lemme \ref{lem:solD1}, il n'y a pas de spineur harmonique
$L^{2}$ pour $|K|>L-4$. L'espace $(S_{4}S_{L})^{S^{1}}$ est alors
de dimension 5, engendré par $(\sigma_{4}=s_{4}Y^{2}g,\,\sigma_{2}=s_{2}Yg,\,\sigma_{0}=s_{0}g,\,\sigma_{-2}=s_{-2}Xg,\,\sigma_{-4}=s_{-4}X^{2}g)$,
et de même on a une base $(\tau_{2}=t_{2}Yg,\,\tau_{0}=t_{0}g,\,\tau_{-2}=t_{-2}Xg)$
de $(S_{2}S_{L})^{S^{1}}$.

Nous calculons à présent les opérateurs $A$, $B$ et $C$.

\begin{claim*}
On a les identités\[
\begin{array}{lclcl}
A\sigma_{4}=3\sigma_{4}, &  & B\sigma_{4}=(K-4)\sigma_{4}, &  & C\sigma_{4}=\frac{L(L+2)-(K-2)(K-4)}{16}(\sigma_{2}-\tau_{2}),\\
A\sigma_{2}=\tau_{2}, & \quad & B\sigma_{2}=\frac{K-2}{2}(\sigma_{2}+\tau_{2}), & \quad & C\sigma_{2}=\sigma_{4}+\frac{L(L+2)-K(K-2)}{16}(\sigma_{0}-\tau_{0}),\\
A\sigma_{0}=-\sigma_{0}, &  & B\sigma_{0}=K\tau_{0}, &  & C\sigma_{0}=\frac{1}{2}(3\sigma_{2}+\tau_{2})+\frac{1}{4}(\sigma_{-2}-\tau_{-2}).\end{array}\]

\end{claim*}
La vérification de ces formules est laissée au lecteur.

Nous pouvons maintenant calculer, pour un spineur $\sigma=\sum_{-4}^{4}a_{i}\sigma_{i}$,
la projection de $\partial_{r}\cdot\mathcal{D}\sigma$ sur $\tau_{2}\in(S_{2}S_{L})^{S^{1}}$.
La formule pour l'opérateur de Dirac indique que seuls les composantes
$\sigma_{0},\sigma_{2},\sigma_{4}$ contribuent et on obtient\begin{equation}
\pi_{\tau_{2}}P\sigma=\left(\frac{u}{2}-\frac{K-2}{2}\right)a_{2}-2\sqrt{1+u}\left(-\frac{L(L+2)-(K-2)(K-4)}{16}a_{4}+\frac{1}{2}a_{0}\right).\label{eq:contrainte}\end{equation}
Cette égalité est une contrainte exprimant $a_{0}$ en fonction de
$a_{2}$ et $a_{4}$, quand $\sigma$ satisfait l'équation $\mathcal{D}\sigma=0$.

Calculons également les projections de $P\sigma$ sur $\sigma_{4}$
et $\sigma_{2}$ : \begin{align}
\pi_{\sigma_{4}}P\sigma & =-2u(1+u)\partial_{u}a_{4}+(-6-\frac{11}{2}u+\frac{3}{2}u-(K-4))a_{4}-2\sqrt{1+u}a_{2},\label{eq:a4}\\
\pi_{\sigma_{2}}P\sigma & =-2u(1+u)\partial_{u}a_{2}+(-6-\frac{11}{2}u-\frac{K-2}{2})a_{2}\label{eq:a2}\\
 & \qquad-2\sqrt{1+u}(\frac{L(L+2)-(K-2)(K-4)}{16}a_{4}+\frac{3}{2}a_{0}).\nonumber \end{align}
Remplaçant $a_{0}$ et $a_{2}$ en fonction de $a_{4}$ grâce aux
équations (\ref{eq:contrainte}) et (\ref{eq:a4}), l'équation (\ref{eq:a2})
devient \[
u(u+1)\partial_{u}^{2}a_{4}+(7u+6)\partial_{u}a_{4}+\frac{-(u+1)L(L+2)+(K^{2}-4K+60)u+24}{4u(u+1)}a_{4}=0\]
dont les deux solutions s'expriment en terme des fonctions hypergéométriques
(voir \cite[chapitre XIV]{WhiWat27}) :\[
u^{-\frac{L}{2}-3}(1+u)^{\frac{K}{2}-1}F(\frac{K-L}{2}-2,\frac{K-L}{2};-L;-u)\,,\]
\[
u^{\frac{L}{2}-2}(1+u)^{\frac{K}{2}-1}F(\frac{K+L}{2}-1,\frac{K+L}{2}+1;L+2;-u)\,.\]
En se rappelant que $u=\sinh^{2}r$, ces formules fournissent un comportement
asymptotique en $0$ en $r^{-L-6}$ et $r^{L-4}$, c'est-à-dire le
comportement trouvé sur les composantes $S_{L+4}$ et $S_{L-4}$ dans
la formule (\ref{eq:cmpt0}). La première des deux solutions ne se
prolonge pas à l'origine, et on en déduit la formule, pour une constante
$A_{4}$,\[
a_{4}(r)=A_{4}\sinh(r)^{L-4}\cosh(r)^{K-2}F(\frac{K+L}{2}-1,\frac{K+L}{2}+1;L+2;-\sinh^{2}r).\]
Sachant que si $a<b$ et $z\rightarrow+\infty$, on a \cite[14.51]{WhiWat27}\[
F(a,b;c;-z)\sim z^{-a}\frac{\Gamma(c)\Gamma(b-a)}{\Gamma(b)\Gamma(c-a)}\,,\]
on déduit que, pour $r\rightarrow+\infty$, \[
a_{4}(r)\sim\frac{1}{L+2}{{L+2 \choose (K+L)/2}}\frac{A_{4}}{\sinh(r)^{4}}\,.\]
On voit là encore se confirmer le comportement asymptotique à l'infini
prédit par le lemme \ref{lem:ds}.

\subsubsection{Conclusion}

Du calcul précédent, on déduit le lemme suivant.

\begin{lem}
\label{lem:injectivite}L'application associant à un spineur $s$,
harmonique $L^{2}$ sur $\bC H^{2}$, sa valeur à l'infini $s_{\infty}$,
est injective.\qed
\end{lem}
Résumons à présent les résultats précédents grâce à la définition
suivante.

\begin{defn}
\label{def:Jplus}Notons $S_{4}^{-}(\pm4)$ les spineurs de poids
$\pm4i$ pour $\sigma_{1}^{-}$. Soit $\rho$ une représentation irréductible
de $U_{1}SU_{2}$, on notera $J_{\rho}^{+}\subset V_{\rho}\otimes(V_{\rho}S_{4}^{-}(\pm4))^{S^{1}}$
l'espace des valeurs à l'infini des spineurs harmoniques dans $\mathcal{S}_{4}^{-}$
provenant de la représentation $\rho$. On notera $\mathcal{J}^{+}$
l'espace des sections de $\mathcal{J}$ sur $S^{3}$ engendré par
les $J_{\rho}^{+}$.
\end{defn}
Les lemmes \ref{lem:solD1} et \ref{lem:injectivite} indiquent que
$J_{\rho}^{+}$ est non nul seulement quand $|K|\leq L-4$, et qu'il
est alors de dimension $\dim J_{\rho}^{+}=\dim V_{\rho}$.

\begin{rem}
\label{rem:continuite}Nous n'avons défini aucune topologie sur $\mathcal{J}^{+}$.
À l'occasion, on précisera toujours la régularité des sections concernées,
par exemple $L^{2}(\mathcal{J}^{+})$ pour les sections $L^{2}$.

Un peu plus de travail après le lemme \ref{lem:injectivite} permet
de montrer qu'en réalité l'application $\sigma_{\infty}\rightarrow\sigma$
qui à une valeur à l'infini dans $\mathcal{J}^{+}$ associe le spineur
harmonique correspondant est continue entre espaces $L^{2}$. Dans
la suite de cet article, nous n'aurons pas besoin de ce fait, les
spineurs harmoniques arrivant a priori avec une régularité suffisante.
\end{rem}

\section{Métriques ACH}

\subsection{Les métriques asymptotiquement complexe-hyperboliques}

Dans cette section, nous reprenons quelques bases sur les métriques
ACH (asymptotiquement complexe-hyperboliques), étudiées dans \cite{Biq00}
et \cite{BiqHer}, avec la perspective de généraliser les énoncés
au cas où la donnée sur le bord n'est plus $C^{\infty}$.

\subsubsection{Métriques ACH}

La métrique hyperbolique complexe peut s'écrire \[
g_{\bC H^{2}}=dr^{2}+\sinh^{2}(r)\gamma_{0}+\sinh^{2}(2r)\eta^{2},\]
où $\eta$ est la forme de contact standard sur la sphère $S^{3}$,
et $\gamma_{0}(\cdot,\cdot)=d\eta(\cdot,J_{0}\cdot)$ est la métrique
dans les directions de contact, provenant de la structure complexe
$J_{0}$ dans les directions de contact.

Rappelons que les déformations $J$ de la structure CR sur les directions
de contact sont paramétrées par les $\phi\in\Omega_{J_{0}}^{0,1}\otimes T_{J_{0}}^{1,0}$,
de sorte que \begin{equation}
T_{J}^{0,1}=\{ X+\phi_{X},X\in T_{J_{0}}^{0,1}\}.\label{eq:phi}\end{equation}
 Un tel $\phi$ est équivalent à la donnée de l'endomorphisme anti-$J_{0}$-linéaire
$\phi+\bar{\phi}$ de $\ker\eta$, qui est une section du fibré $\mathcal{J}$
des endomorphismes symétriques à trace nulle de $\ker\eta$.

À une structure CR $J$ est associée une métrique $\gamma$ dans les
directions de contact. On peut alors construire une métrique ayant
asymptotiquement le même comportement que la métrique hyperbolique
complexe, à savoir \begin{equation}
g_{0}=dr^{2}+e^{2r}\gamma+e^{4r}\eta^{2}.\label{eq:g0}\end{equation}
Plus généralement \cite{Biq00}, on dira que la métrique $g$ est
ACH (asymptotiquement complexe-hyperbolique) si, près de l'infini,
\begin{equation}
g=dr^{2}+e^{2r}\gamma+e^{4r}\eta^{2}+O(e^{-r}),\label{eq:ACH}\end{equation}
où le $O(e^{-r})$ est mesuré par rapport à la métrique $g_{0}$ (on
demande aussi que les dérivées satisfassent la même décroissance).
La courbure de $g$ à l'infini s'approche alors à l'ordre $O(e^{-r})$
de celle de $g_{0}$, et $J$ (ou $\gamma$) est appelé l'infini conforme
de $g$.

\begin{example}
Voyons la métrique hyperbolique complexe comme la métrique de Bergmann
sur la boule dans $\bC^{2}$. Si on déforme la boule en un domaine
pseudoconvexe de $\bC^{2}$, on peut construire une métrique Kähler-Einstein,
la métrique de Cheng-Yau \cite{CheYau80}, qui est ACH.
\end{example}

\subsubsection{Métriques d'Einstein\label{sec:met-Einstein}}

Plus généralement, toute petite déformation $J$ de la structure CR
$J_{0}$ est l'infini conforme d'une métrique ACH Einstein \cite[chapitre I]{Biq00}. 

Le comportement asymptotique à un ordre élevé de cette métrique d'Einstein
a été étudié dans \cite[théorème 3.3 et corollaire 3.4]{BiqHer} :
plus précisément, une structure CR $J$ est toujours le bord d'une
structure complexe formelle locale ; étant donné $J$ sur le bord,
on peut construire un développement explicite d'ordre 4 pour une métrique
ACH $\bar{g}$, d'infini conforme $J$, Kähler-Einstein à l'ordre
$4+\delta$ près de l'infini ($0<\delta<1$). Cette métrique a la
forme suivante : le terme principal est $g_{0}$ défini par (\ref{eq:g0}),
puis on obtient \begin{equation}
\bar{g}=g_{0}+e^{-2r}g_{2}+e^{-3r}g_{3}+e^{-4r}g_{4}+O(e^{-5r}),\label{eq:gbar}\end{equation}
où les termes $g_{i}$ sont des termes sur le bord, obtenus non linéairement
à partir de $J$, mais impliquant au plus $i$ dérivées horizontales
de $J$. Il s'agit donc d'un développement partiel \og polyhomogène \fg{}.
Cette métrique est d'Einstein à un ordre élevé, $\Ric^{\bar{g}}+6\bar{g}=O(e^{-5r})$,
et le terme d'ordre $e^{-5r}$ implique au plus $5$ dérivées horizontales
de $J$ ($\Ric^{\bar{g}}$ lui-même implique $6$ dérivées de $J$,
mais la dérivée sixième n'apparaît qu'à l'ordre $6$).

La métrique asymptotiquement Kähler-Einstein $\bar{g}$ est une excellente
approximation de la métrique d'Einstein $g$ qui remplit $J$. Plus
précisément \cite[corollaire 5.4]{BiqHer}, \begin{equation}
g=\overline{g}+\Gamma e^{-2r}+O(e^{-(4+\delta)r}),\label{eq:ggbar}\end{equation}
où $\Gamma$ est une section sur $S^{3}$ du fibré $\mathcal{J}$
des endomorphismes symétriques à trace nulle de $\ker\eta$ (donc
$\Gamma e^{-2r}$ est un terme homogène d'ordre 4) ; si $J$ est $C^{\infty}$,
alors il en est de même pour $\Gamma$. 

Dans la suite de cette section, nous allons développer quelques outils
pour comprendre la régularité de la métrique $g$ quand $J$ n'est
plus lisse.

\subsubsection{Espaces fonctionnels\label{sec:fonctionnels}}

Commençons par introduire quelques espaces fonctionnels sur le bord.
On y dispose de la connexion de Webster, $\nabla^{W}$. Celle-ci peut
être étendue à l'intérieur de la manière suivante : pour une métrique
ACH, ayant le comportement asymptotique (\ref{eq:ACH}), la dérivée
covariante $\nabla_{\partial_{r}}$ permet d'identifier les fibrés
tangents sur les sphères concentriques (les vecteurs parallèles sont
approximativement $\partial_{r}$, $e^{-2r}R$, et $e^{-r}h$ pour
$h$ un vecteur de contact) ; le tiré en arrière de $\nabla^{W}$,
que l'on notera à nouveau $\nabla^{W}$, définit alors une connexion
à l'intérieur. Par \cite[Corollaire 2.3]{BiqHer}, on a en réalité
\begin{equation}
\nabla=\nabla^{W}+a_{0}+O(e^{-r}),\quad\nabla^{W}a_{0}=0;\label{eq:a0}\end{equation}
dans cette formule, les dérivées du $O(e^{-r})$ ont la même décroissance.
Dans le cas particulier de $\bC H^{2}$, on a la formule explicite
(\ref{eq:nabla1}) dans laquelle les termes $\rho(Y_{i})$ fournissent
la connexion de Webster et les termes $\rho_{E}(Y_{i})$ fournissent
la correction $a_{0}$ ; cette correction $a_{0}$ est la même pour
toutes les métriques ACH.

On notera $\square$ le laplacien (hypoelliptique) associé à la dérivation
dans les directions horizontales,\[
\square f=-\sum_{i=1}^{2}(\nabla^{W})_{h_{i},h_{i}}^{2}f,\]
où $(h_{1},h_{2})$ est une base orthonormale du plan de contact.
Enfin, on notera $FS^{k}$ l'espace des fonctions (ou sections d'un
fibré) sur $S^{3}$, ayant $k$ dérivées horizontales dans $L^{2}$
: ce sont les espaces de Folland-Stein \cite{FolSte74}. Plus généralement,
nous aurons besoin d'espaces avec un nombre de dérivées fractionnaire,
$FS^{\delta}$, que l'on peut définir à partir de la théorie spectrale
de $\square$ : notons $\lambda^{2}$ ses valeurs propres ($\lambda\geq0$),
si $f=\sum f_{\lambda}$ est la décomposition de $f$ sur les espaces
propres correspondants de $\square$, alors l'espace $FS^{\delta}$
est défini par la norme $\sum_{\lambda}(1+\lambda)^{2\delta}\| f_{\lambda}\|^{2}$.

Passons maintenant aux espaces fonctionnels sur $\bC H^{2}$. Nous
avons l'espace de Sobolev classique $H^{k}$ des fonctions ayant $k$
dérivées dans $L^{2}$. Nous utiliserons aussi la version à poids,
en notant $H_{\delta}^{k}$ l'espace des fonctions $f$ telles que
$e^{(\delta-2)r}f\in H^{k}$ ; la convention est choisie de sorte
que la fonction $e^{-dr}$ soit dans $H_{\delta}^{k}$ si et seulement
si $d>\delta$.

Dans ces espaces de Sobolev, la dérivation suivant une direction horizontale
de $S^{3}$ est toujours affectée d'un poids $e^{-r}$. Or il sera
utile de considérer des espaces fonctionnels avec une meilleure régularité
horizontale. Cela nous mène à définir, pour $\ell\leq k$, l'espace
$H_{\delta}^{k;\ell}$ des fonctions $f\in H_{\delta}^{k}$ telles
que $\chi(r)(\nabla_{h_{i}}^{W})^{\ell}f\in H_{\delta}^{k-\ell}$,
où $\chi(r)$ est une fonction de coupure, telle que $\chi(r)=1$
pour $r\geq R+1$ et $\chi(r)=0$ pour $r\leq R$, pour un $R>0$
suffisamment grand.

Un point important est qu'on peut supposer que les diverses structures
CR sur le bord (proches de la structure standard) ont la même structure
de contact sous-jacente : il en résulte que les espaces fonctionnels
définis plus haut ne dépendent pas de la métrique ACH choisie. En
particulier, on peut utiliser la métrique hyperbolique complexe $g_{\bC H^{2}}$
et la structure CR standard $J_{0}$ sur le bord pour les définir.

Dans toute la suite, \emph{on supposera $\ell$ fixé très grand},
de sorte que $FS^{\ell}\subset C^{0}$ ; en particulier, cela implique
que les espaces $FS^{\ell}$,...,$FS^{\ell+4}$ que nous considérerons
soient stables par multiplication, c'est-à-dire soient des algèbres.
Nous demandons la même chose pour les espaces de Sobolev $H^{k}$.

\begin{rem}
\label{rem:dR}Si on a une fonction $f\in H_{\delta}^{2;1}$, alors
$Rf=-[h,Jh]f+[h,Jh]_{H}f\in L_{\delta-1}^{2}$, ce qui est mieux que
$e^{-2r}Rf\in L_{\delta}^{2}$. Cela reste valable pour une section
d'un fibré. 
\end{rem}

\subsubsection{Opérateur de prolongement}

Si une fonction $f$ est donnée sur le bord à l'infini, on peut chercher
à prolonger $f$ à l'intérieur en une fonction $E_{0}f$ ayant la
meilleure régularité possible. Cela est réalisé à travers la construction
suivante. 

Partons de la décomposition de $f$ sur les espaces propres de $\square$,
et posons \[
E_{0}f=\chi(r)\sum_{\lambda}f_{\lambda}e^{-\lambda e^{-r}},\]
où, comme plus haut, $\chi$ est une fonction de coupure valant $1$
près de l'infini et $0$ près de $0$.

\begin{lem}
\label{lem:reg0}Si $f\in FS^{\delta}$ avec $0<\delta<1$, alors
$f-E_{0}f\in L_{\delta}^{2}$ et $\nabla(E_{0}f)\in H_{\delta}^{\infty}$
(avec normes contrôlées par $\| f\|_{FS^{\delta}}$). 
\end{lem}
\begin{rem}
Il est clair que $[E_{0},\square]=0$. Par conséquent, si $f\in FS^{\ell+\delta}$,
alors $\square^{\frac{\ell}{2}}E_{0}f-\square^{\frac{\ell}{2}}f\in L_{\delta}^{2}$.
\end{rem}
\begin{proof}
Posons $g=\partial_{r}\sum f_{\lambda}e^{-\lambda e^{-r}}=\sum\lambda e^{-r}f_{\lambda}e^{-\lambda e^{-r}}$,
alors \[
\| g\|_{L_{\delta}^{2}(r\geq R)}^{2}=\sum_{\lambda}\int_{r\geq R}\lambda^{2}e^{-2r}|f_{\lambda}|^{2}e^{-2\lambda e^{-2r}}e^{2(-2+\delta)r}\textrm{vol}^{g}.\]
En utilisant le fait que $\textrm{vol}^{g}\sim e^{4r}$, on est ramené
à contrôler, en utilisant $\delta<1$, \begin{align}
\sum_{\lambda}\| f_{\lambda}\|^{2}\int_{R}^{+\infty}\lambda^{2}e^{2(-1+\delta)r}e^{-2\lambda e^{-2r}}dr & =\sum_{\lambda}\lambda^{2\delta}\| f_{\lambda}\|^{2}\int_{0}^{\lambda e^{-R}}u^{2(1-\delta)}e^{-2u}\frac{du}{u}\label{eq:int-lambda}\\
 & \leq c\sum_{\lambda}\lambda^{2\delta}\| f_{\lambda}\|^{2}\nonumber \end{align}
d'où résulte l'estimation \[
\| g\|_{L_{\delta}^{2}(r\geq R)}^{2}\leq c\| f\|_{FS^{\delta}}^{2},\]
qui implique en particulier $f-E_{0}f\in L_{\delta}^{2}$. 

Une dérivation tangentielle dans une direction de contact (toujours
accompagnée d'un facteur $e^{-r}$ dans $\nabla E_{0}f$) se traduit
par la multiplication du terme $f_{\lambda}e^{-\lambda e^{-r}}$ par
$\lambda e^{-r}$, et donc satisfait la même estimation. Plus généralement,
une dérivation à l'ordre $k$ dans les directions de contact équivaut
à une multiplication par $\lambda^{k}e^{-kr}$, se traduisant par
le remplacement de $u^{2(1-\delta)}$ par $u^{2(k-\delta)}$ dans
l'intégrale (\ref{eq:int-lambda}) : on obtient ainsi un contrôle
similaire. Le lemme en découle.
\end{proof}
Il y a une réciproque qui montre que les valeurs au bord de régularité
$FS^{\delta}$ correspondent exactement aux fonctions dont les dérivées
décroissent en $e^{-\delta r}$ :

\begin{lem}
\label{lem:rpq}Si une fonction $f$ satisfait $df\in L_{\delta}^{2}$
alors $f$ a une valeur à l'infini $f_{\infty}\in FS^{\delta}$, de
sorte que $f-f_{\infty}\in L_{\delta}^{2}$.
\end{lem}
\begin{proof}
L'existence d'une valeur à l'infini $f_{\infty}$ telle que $f-f_{\infty}\in L_{\delta}^{2}$
est immédiate. Prouvons l'estimation $FS^{\delta}$ : soit $\chi(r)$
une fonction de coupure, telle que $\chi(r)=1$ pour $r\geq R+1$
et $\chi(r)=0$ pour $r\leq R$. Alors \[
f_{\infty}=-\int_{R}^{+\infty}\partial_{r}(\chi f)=-\int_{R}^{+\infty}(\partial_{r}\chi)f+\chi\partial_{r}f\]
d'où l'estimation \[
\frac{1}{4}|f_{\infty}|^{2}\leq e^{2(1-\delta)R}\int f^{2}e^{2(\delta-1)r}+e^{-2\delta R}\int|\partial_{r}f|^{2}e^{2\delta r}.\]
Restreignons à l'espace propre de $\square$ pour la valeur propre
$\lambda^{2}$, et choisissons $R=\ln\lambda$, alors l'estimation
devient \[
\frac{1}{4}\| f_{\infty}\|^{2}\leq\lambda^{-2\delta}\left(\| e^{-r}\lambda f\|_{L_{\delta}^{2}}^{2}+\|\partial_{r}f\|_{L_{\delta}^{2}}^{2}\right)\]
qui est exactement l'estimation $FS^{\delta}$ voulue.
\end{proof}
\begin{rem}
Le lemme reste clairement valable si $f$ est une section d'un fibré,
à condition de remplacer la condition $df\in L_{\delta}^{2}$ par
$\nabla^{W}f\in L_{\delta}^{2}$.
\end{rem}
Le lemme \ref{lem:reg0} donne une approximation d'une valeur au bord
$f$ à un ordre $\delta$, mais il peut être utile d'en construire
à un ordre plus élevé : pour obtenir une approximation à l'ordre $j+\delta$,
il suffit d'enlever les termes d'ordre $1$ à $j$ dans $e^{-\lambda e^{-r}}$,
en définissant \[
E_{j}f=\chi(r)\sum_{\lambda}f_{\lambda}(e^{-\lambda e^{-r}}+\lambda e^{-r}-\frac{\lambda^{2}}{2}e^{-2r}-\cdots-\frac{(-\lambda)^{j}}{j!}e^{-jr}).\]
On obtient alors, de manière similaire :

\begin{lem}
\label{lem:reg}Si $f\in FS^{\ell+\delta}$, alors 
\begin{enumerate}
\item $E_{\ell}f-f\in L_{\ell+\delta}^{2}$ ;
\item plus généralement, pour $j\leq\ell$, on a $\square^{\frac{j}{2}}E_{\ell}f-\square^{\frac{j}{2}}f\in L_{\ell-j+\delta}^{2}$
;
\item $\nabla^{\ell+1}E_{\ell}f\in H_{\ell+\delta}^{\infty}$.
\end{enumerate}
\end{lem}
\begin{proof}
La démonstration est la même que celle du lemme \ref{lem:reg0}.
\end{proof}

\subsubsection{Opérateurs elliptiques pour les métriques ACH\label{sec:elliptique}}

Rappelons maintenant brièvement quelques faits sur les opérateurs
elliptiques sur $\bC H^{2}$, ou plus généralement pour une métrique
ACH. Ces opérateurs sont redevables de la théorie des \og edge operators
 \fg{}, voir \cite{EpsMelMen91} dans le cas scalaire pour la métrique
hyperbolique complexe. Soit un opérateur homogène \[
P=\nabla^{*}\nabla+\mathcal{R},\]
où $\mathcal{R}$ est un terme d'ordre 0, supposé inversible dans
$L^{2}$, c'est-à-dire satisfaisant une estimation $(Pu,u)\geq c\| u\|^{2}$.
Le comportement de $P$ est alors gouverné par \og l'opérateur indiciel
 \fg{}, constitué des termes d'ordre 0 à l'infini : c'est un opérateur
du type \[
-\partial_{r}^{2}-4\partial_{r}+A,\]
où $A$ est un opérateur d'ordre 0. Soit $\lambda$ la plus petite
valeur propre de $A$, alors les poids critiques de l'opérateur indiciel
(c'est-à-dire les poids $\delta$ tels que $\exp(-\delta r)$ est
dans le noyau de l'opérateur indiciel) sont $\delta_{\pm}=2\pm\sqrt{4+\lambda}$,
et $P$ est un isomorphisme \begin{equation}
P:H_{\delta}^{k+2}\stackrel{\sim}{\longrightarrow}H_{\delta}^{k}\label{eq:isoHk}\end{equation}
pour tout poids $\delta$ compris entre les poids critiques, $\delta_{-}<\delta<\delta_{+}$,
voir \cite[proposition I.2.5]{Biq00} (dans cette proposition, l'isomorphisme
est montré dans les espaces de Hölder, à partir d'une estimation $L^{2}$
(lemme I.2.3) qui implique aussi immédiatement l'isomorphisme dans
les espaces de Sobolev) ; ce résultat est valable aussi bien pour
la métrique hyperbolique complexe $g_{\bC H^{2}}$ que plus généralement
pour une métrique ACH \cite[I.3]{Biq00}.

Si on résout le problème $Pu=v$ avec $v$ décroissant plus rapidement
que $e^{-\delta_{+}r}$, alors la forme de l'opérateur indiciel suggère
que \[
u=u_{\infty}e^{-\delta_{+}r}+u',\]
avec $u_{\infty}$ section de l'espace propre $\mathcal{U}$ de $A$
pour la valeur propre $\lambda$, et $u'$ décroissant plus rapidement
que $e^{-\delta_{+}r}$ : intrinsèquement, on obtient $u_{\infty}$
par la limite \[
u_{\infty}=\lim_{r\rightarrow+\infty}e^{\delta_{+}r}u.\]
La régularité de la limite est établie dans le théorème suivant :

\begin{thm}
\label{thm:regularite}Soit $0<\delta<1$. Supposons que l'opérateur
indiciel $P$ sur $\bC H^{2}$ n'ait pas d'autre poids critique entre
$\delta_{+}$ et $\delta_{+}+\delta$.

1) Si $v\in H_{\delta_{+}+\delta}^{k}$ alors $u_{\infty}\in FS^{\delta}(\mathcal{U})$
et $u-e^{-\delta_{+}r}E_{0}u_{\infty}\in H_{\delta_{+}+\delta}^{k+2}$. 

2) Plus généralement, on obtient ainsi un isomorphisme\[
P:FS^{\ell+\delta}(\mathcal{U})\oplus H_{\delta_{+}+\delta}^{k+2;\ell}\rightarrow H_{\delta_{+}+\delta}^{k;\ell},\]
 défini par $(u_{\infty},u)\rightarrow P(e^{-\delta_{+}r}E_{0}u_{\infty}+u)$.
\end{thm}
C'est le type de résultat que l'on s'attend à voir dériver de la théorie
des \og edge operators \fg{} : dans le cas hyperbolique réel, cela
est établi dans \cite[théorème 7.14]{Maz91}, et dans le cas hyperbolique
complexe qui nous intéresse, le calcul pseudodifférentiel a été développé
dans \cite{EpsMelMen91}. Mais il semble difficile de trouver une
référence qui ne se limite pas aux problèmes scalaires, auxquels ne
peuvent a priori pas se réduire les systèmes que nous considérons.
Pour être complet, nous incluons ici une démonstration élémentaire,
similaire à celle donnée dans \cite[proposition 5.2]{BiqHer} dans
le cas $C^{\infty}$, valable sous les hypothèses techniques suivantes
:

\begin{enumerate}
\item $k>\ell$ : cette limitation est sans importance pour nous, car les
solutions d'équations elliptiques que nous regardons sont toujours
localement $C^{\infty}$, de sorte que l'on peut prendre $k$ aussi
grand que l'on voudra ;
\item $\delta_{+}-\delta_{-}>2-\delta$ : dans les cas que nous regardons,
le couple $(\delta_{-},\delta_{+})$ est égal à celui du laplacien
scalaire, soit $(0,4)$.
\end{enumerate}
\begin{rem}
Le théorème est énoncé uniquement dans le cas de $\bC H^{2}$, le
seul utilisé dans cet article : cela permet de simplifier considérablement
la démonstration. Cependant il a une portée plus générale : si l'opérateur
$P$ est un isomorphisme dans $L^{2}$ pour une certaine métrique
ACH, alors le théorème demeure valable, à condition d'avoir pour les
coefficients de la métrique une régularité suffisante (dépendant de
la valeur de $\ell$). Plus généralement, si l'opérateur $P$ pour
une métrique ACH n'est que Fredholm entre les poids $\delta_{-}$
et $\delta_{+}$, alors le théorème \ref{thm:regularite} doit être
remplacé par l'assertion que $P$ reste Fredholm, avec indice inchangé,
entre les espaces $FS^{\ell+\delta}(\mathcal{U})\oplus H_{\delta_{+}+\delta}^{k+2;\ell}\rightarrow H_{\delta_{+}+\delta}^{k;\ell}$
(la démonstration à partir du résultat dans $\bC H^{2}$ se fait comme
dans \cite[I.3.B]{Biq00}).
\end{rem}
\begin{proof}
Tout d'abord, le champ de vecteurs $R$ est une isométrie infinitésimale
de $\bC H^{2}$ ; comme $P$ est homogène, on en déduit \[
[P,\mathcal{L}_{R}]=[P,\nabla_{R}^{W}]=0.\]
D'autre part, on a la commutation suivante \cite[lemme 2.6]{BiqHer}
: soit $h$ un champ de vecteurs (de norme $1$) dans la distribution
de contact, alors\begin{equation}
[P,\nabla_{h}^{W}]=2\nabla_{e^{-r}Jh}\nabla_{e^{-r}R}^{W}+Q,\label{eq:com2}\end{equation}
où $Q$ est un opérateur différentiel d'ordre 2, préservant les poids.
Plus précisément, $Q$ est élément d'une algèbre d'opérateurs différentiels
$\mathcal{Q}$, telle que :
\begin{itemize}
\item l'algèbre $\mathcal{Q}$ contient les dérivations horizontales $\nabla_{e^{-r}h}^{W}$
;
\item si $Q\in\mathcal{Q}$ est d'ordre $j$, alors il est continu $H_{\delta'}^{k}\rightarrow H_{\delta'}^{k-j}$
(pour tout $\delta'$) et $FS^{\ell+\delta}\oplus H_{\delta_{+}+\delta}^{k;\ell}\rightarrow H_{\delta_{+}+\delta}^{k-j;\ell}$
(concrètement, les termes de $Q$, d'homogénéïté $0$ par rapport
aux poids, contiennent des dérivations horizontales, envoyant ainsi
un terme $e^{-\delta_{+}r}E_{0}u_{\infty}$ dans un espace à décroissance
en $e^{-(\delta_{+}+\delta)r}$).
\end{itemize}
Supposons donc l'opérateur différentiel $P$ inversible entre les
poids $\delta_{-}$ et $\delta_{+}$. Soit $\delta>0$, $v\in H_{\delta_{+}+\delta}^{1}$
et $u$ la solution de $Pu=v$, de sorte que $u\in H_{\delta_{+}-\varepsilon}^{3}$
pour tout $\varepsilon>0$ petit. Tout d'abord, \[
P\nabla_{R}^{W}u=\nabla_{R}^{W}Pu\in L_{\delta_{+}-2+\delta}^{2}.\]
Utilisant une fonction de coupure $\chi(r)$ nulle en dehors de $r\geq R$,
on déduit $P\chi\nabla_{R}^{W}u\in L_{\delta_{+}-2+\delta}^{2}$ et
par conséquent (utilisant la condition $\delta_{+}-2+\delta>\delta_{-}$)
$\chi\nabla_{R}^{W}u\in H_{\delta_{+}-2+\delta}^{2}$, d'où $\nabla_{e^{-2r}R}^{W}u\in H_{\delta_{+}+\delta}^{2}$,
ce qui est le contrôle attendu sur cette dérivée.

Pour simplifier l'écriture de la suite du raisonnement, on omettra
à chaque fois l'utilisation de la fonction $\chi$.

On a, d'après la commutation (\ref{eq:com2}), \[
P\nabla_{h}^{W}u=[P,\nabla_{h}^{W}]u+\nabla_{h}^{W}Pu\in L_{\delta_{+}-1+\delta}^{2}.\]
Il en résulte $\nabla_{h}^{W}u\in H_{\delta_{+}-1+\delta}^{2}$, soit
$\nabla_{e^{-r}h}^{W}u\in H_{\delta_{+}+\delta}^{2}$, ce qui est
à nouveau le contrôle attendu. La forme de l'opérateur $P=\nabla^{*}\nabla+\mathcal{R}$
nous indique alors que ($A$ désignant l'opérateur linéaire formé
des termes d'ordre $0$ à l'infini de $P$) \[
Pu-(-\partial_{r}^{2}u-4\partial_{r}u+Au)\in L_{\delta_{+}+\delta}^{2}\]
 et donc \[
-\partial_{r}^{2}u-4\partial_{r}u+Au\in L_{\delta_{+}+\delta}^{2}.\]
Des arguments élémentaires sur les équations différentielles ordinaires
impliquent la convergence à l'infini de $e^{\delta_{+}r}u$ vers une
limite $u_{\infty}$ dans le sous-espace propre de $A$ correspondant
à ce poids, et \[
\partial_{r}(u-e^{-\delta_{+}r}u_{\infty})\in L_{\delta_{+}+\delta}^{2}.\]
De ce contrôle et des contrôles précédents sur $\nabla^{W}u$ résulte
finalement $\nabla^{W}(e^{\delta_{+}r}u)\in L_{\delta}^{2}$, d'où
la régularité $u_{\infty}\in FS^{\delta}$ par le lemme \ref{lem:rpq}.
Le lemme \ref{lem:reg0} indique alors que $P(e^{-\delta_{+}r}E_{0}u_{\infty})\in H_{\delta_{+}+\delta}^{\infty}$.

Finalement, on déduit que $P(u-e^{-\delta_{+}r}E_{0}u_{\infty})\in H_{\delta_{+}+\delta}^{1}$
et $u-e^{-\delta_{+}r}u_{\infty}\in L_{\delta_{+}+\delta}^{2}$. La
régularité elliptique locale pour l'opérateur $P$ donne alors $u-e^{-\delta_{+}r}E_{0}u_{\infty}\in H_{\delta_{+}+\delta}^{3}$.
Si de plus $v\in H_{\delta_{+}+\delta}^{k}$, alors on obtient $P(u-e^{-\delta_{+}r}E_{0}u_{\infty})\in H_{\delta_{+}+\delta}^{k}$
d'où résulte $u-e^{-\delta_{+}r}E_{0}u_{\infty}\in H_{\delta_{+}+\delta}^{k+2}$.
Cela démontre le théorème dans le cas $\ell=0$.

Supposons à présent que l'on ait la régularité supplémentaire $v\in H_{\delta_{+}+\delta}^{k;\ell}$
avec $k>\ell$. Commençons par $\ell=1$. D'après la remarque \ref{rem:dR},
on a $\nabla_{R}^{W}v\in H_{\delta_{+}+\delta-1}^{k-2}$ et donc $P\nabla_{R}^{W}u\in H_{\delta_{+}+\delta-1}^{k-2}$,
d'où résulte $\nabla_{R}^{W}u\in H_{\delta_{+}+\delta-1}^{k}$. Appliquons
ensuite la commutation (\ref{eq:com2}) : \begin{equation}
P\nabla_{h}^{W}u=\nabla_{h}^{W}Pu+2\nabla_{e^{-r}Jh}\nabla_{e^{-r}R}^{W}u+Qu\in H_{\delta_{+}+\delta}^{k-1},\label{eq:Pnhu}\end{equation}
d'où résulte $\nabla_{h}^{W}u_{\infty}\in FS^{\delta}$ et donc $u_{\infty}\in FS^{\delta+1}$.
On obtient le cas général par une récurrence sur $\ell$ comme dans
la démonstration de \cite[proposition 5.2]{BiqHer}.
\end{proof}

\subsubsection{Régularité des métriques d'Einstein d'infini conforme $FS^{\ell+4+\delta}$}

Pour étudier le problème d'autodualité pour les métriques ACH Einstein,
nous avons besoin d'étendre l'étude décrite section \ref{sec:met-Einstein}
au cas où l'infini conforme n'est plus $C^{\infty}$, mais plutôt
dans un espace de Folland-Stein $FS^{\ell+4+\delta}$. On prendra
cependant $\ell$ très grand, comme convenu section \ref{sec:fonctionnels}.
En revanche, on pourra garder à l'intérieur des métriques de régularité
locale $C^{\infty}$.

Reprenons donc la construction de la section \ref{sec:met-Einstein},
en supposant que la structure CR sur le bord, $J$, ou de manière
équivalente $g_{0}$, est maintenant de régularité $FS^{\ell+4+\delta}$.
Grâce à la régularisation des données au bord par les lemmes \ref{lem:reg0}
et \ref{lem:reg}, on peut construire la métrique $\bar{g}$, qui
est Kähler-Einstein à un ordre élevé, en modifiant (\ref{eq:gbar})
par \begin{equation}
\overline{g}=E_{4}g_{0}+E_{2}g_{2}e^{-2r}+E_{1}g_{3}e^{-3r}+E_{0}g_{4}e^{-4r}.\label{eq:gbar2}\end{equation}
C'est une métrique qui coïncide avec (\ref{eq:gbar}) à l'ordre $4+\delta$
; le tenseur de Ricci a la régularité attendue :

\begin{lem}
\label{lem:Ricgbar2}Pour la métrique $\overline{g}$ définie par
(\ref{eq:gbar2}), le tenseur de Ricci satisfait

\[
\Ric^{\overline{g}}+6\overline{g}\in H_{4+\delta}^{\infty;\ell}.\]
Plus généralement, la courbure de $\overline{g}$ a le développement
\[
R^{\overline{g}}=R_{0}+\cdots+e^{-4r}R_{4}+O(e^{-(4+\delta)r}),\]
avec $R_{i}\in FS^{\ell+4-i+\delta}$ et le $O(e^{-(4+\delta)r})$
signifie que $R^{\overline{g}}-(E_{4}R_{0}+\cdots+e^{-4r}E_{0}R_{4})\in H_{4+\delta}^{\infty;\ell}.$
\end{lem}
\begin{proof}
Tout d'abord, l'assertion sur $\Ric^{\overline{g}}$ est une conséquence
de l'assertion sur $R^{\overline{g}}$, puisque $\overline{g}$ a
été construite de sorte que les termes du développement polyhomogène
de $\Ric^{\overline{g}}+6\overline{g}$ s'annulent jusqu'à l'ordre
$4$.

Montrons donc l'assertion sur $R^{\overline{g}}$. Si $g_{0}$ est
lisse, le résultat est clair. En revanche, si $g_{0}$ est seulement
$FS^{\ell+4+\delta}$, il n'est pas immédiat de voir comment les régularisations
$E_{i}$ passent à travers la courbure. Ramenons-nous à un problème
linéaire : si on a un chemin de métriques au bord $g_{0}(t)$, tel
que $g_{0}(0)$ est lisse et $g_{0}(1)=g_{0}$, il suffit de montrer
que $\frac{d}{dt}R^{\overline{g}(t)}$ a aussi un développement polyhomogène
du type voulu. Or, \[
d_{g}R(\dot{h})=F(\nabla_{g}^{2}\dot{h},R_{g}\odot\dot{h}),\]
où $F$ est un opérateur linéaire et $\odot$ est une opération bilinéaire.
Ayant pris $\ell$ assez grand, les développements polyhomogènes souhaités
sont stables par multiplication, d'où on déduit que, si $R_{g}$ et
$\dot{h}$ ont un développement polyhomogène comme dans le lemme,
alors $d_{g}R(\dot{h})$ admet un développement du même type.
\end{proof}
Si $J$ est assez proche de $J_{0}$, c'est l'infini conforme d'une
métrique ACH, Einstein (section \ref{sec:met-Einstein}). Nous raffinons
cette étude pour obtenir :

\begin{lem}
\label{lem:regg}Si $J$, de régularité $FS^{\ell+4+\delta}$ sur
$S^{3}$, est suffisamment proche de $J_{0}$, alors il existe une
métrique ACH, Einstein, $g$, d'infini conforme $J$, telle que \[
g=\bar{g}+\Gamma e^{-2r}+G\]
avec $G\in H_{4+\delta}^{\infty;\ell}$, et $\Gamma$ une section
sur $S^{3}$ du fibré $\mathcal{J}$ des endomorphismes symétriques
à trace nulle de $\ker\eta$, de régularité $FS^{\ell+\delta}$ (donc
$\Gamma e^{-2r}$ est un terme d'ordre 4).
\end{lem}
\begin{proof}
Rappelons brièvement la construction de la métrique d'Einstein \cite[chapitre I]{Biq00}
: on résout le problème \[
\Ric^{\overline{g}+h}+6(\overline{g}+h)=0,\quad\delta^{\overline{g}}h+\frac{1}{2}d\tr^{\overline{g}}h=0;\]
la seconde équation fixe en $\overline{g}$ une jauge à l'action des
difféomorphismes induisant l'identité sur le bord à l'infini. Ces
deux équations sont équivalentes, pour $\overline{g}+h$ proche de
la métrique hyperbolique complexe $g_{\bC H^{2}}$, à \[
\Phi^{\overline{g}}(h)=\Ric^{\overline{g}+h}+6(\overline{g}+h)+\delta^{\overline{g}+h}\left(\delta^{\overline{g}}h+\frac{1}{2}d\tr^{\overline{g}}h\right)=0,\]
dont la linéarisation en $g_{\bC H^{2}}$ est, en notant $\Rrond\dot{h}_{X,Y}=\sum\dot{h}(R_{X,e_{i}}Y,e_{i})$
l'action de la courbure de $g_{\bC H^{2}}$ sur les formes quadratiques,\[
d\Phi^{g_{\bC H^{2}}}(\dot{h})=\frac{1}{2}\nabla^{*}\nabla\dot{h}-\Rrond\dot{h.}\]

Cet opérateur homogène sur $\bC H^{2}$ est inversible dans $L^{2}$
\cite[lemme I.1.5]{Biq00}, et relève de la théorie rappelée section
\ref{sec:elliptique}. Son comportement à l'infini est gouverné par
l'opérateur indiciel, ici égal à \[
-\partial_{r}^{2}-4\partial_{r}+A,\]
où $A$ est un opérateur positif, de noyau égal au fibré $\mathcal{J}$
\cite[I.4.B]{Biq00} ; le premier poids critique est donc égal à $4$.
D'après le théorème \ref{thm:regularite}, on a donc un isomorphisme
\[
d\Phi^{g_{\bC H^{2}}}:FS^{\ell+\delta}(\mathcal{J})\oplus H_{4+\delta}^{k+2;\ell}\rightarrow H_{4+\delta}^{k;\ell}.\]

Par le lemme \ref{lem:Ricgbar2}, on a \[
\Phi^{\overline{g}}(0)=\Ric^{\overline{g}+h}+6(\overline{g}+h)\in H_{4+\delta}^{\infty;\ell}\]
donc on peut appliquer le théorème des fonctions implicites à l'équation
$\Phi^{\overline{g}}(h)=0$. Le lemme en résulte.
\end{proof}

\subsection{Noyaux de l'opérateur de Dirac en famille}

Dans la section \ref{sec:1}, on a vu que sur $\bC H^{2}$ le noyau
de l'opérateur de Dirac $\ker\mathcal{D}\subset L^{2}(\mathcal{S}_{-}^{4})$
est de dimension infinie. On va maintenant montrer que tel demeure
le cas quand on perturbe la métrique, et même que ces noyaux forment
un fibré au-dessus de l'espace des métriques ACH.

\subsubsection{L'opérateur de Dirac sur $\mathcal{S}^{+}\mathcal{S}_{3}^{-}$}

L'espace des 1-formes à valeurs dans $\Omega^{-}$ se décompose en
\[
\Omega^{1}\Omega^{-}=\mathcal{S}^{+}\mathcal{S}^{-}\mathcal{S}_{2}^{-}=\mathcal{S}^{+}\mathcal{S}^{-}\oplus\mathcal{S}^{+}\mathcal{S}_{3}^{-}.\]
 Il est classique que l'opérateur de Dirac $\mathcal{D}$ sur $\mathcal{S}^{+}\mathcal{S}_{3}^{-}$
s'identifie alors à $d^{*}+\sqrt{2}d^{-}$, à valeurs a priori dans
$\Omega^{-}+\Omega^{-}\Omega^{-}$, mais en réalité la restriction
à $\mathcal{S}^{+}\mathcal{S}_{3}^{-}$ est à valeurs dans \[
\Omega^{-}+\fS_{0}^{2}\Omega^{-}=\mathcal{S}_{2}^{-}+\mathcal{S}_{4}^{-}=\mathcal{S}^{-}\mathcal{S}_{3}^{-}.\]
 De plus, on a le résultat suivant.

\begin{lem}
\label{lem:W}Soit $\scal$ la courbure scalaire de $\bC H^{2}$,
alors on a sur $\mathcal{S}^{+}\mathcal{S}_{3}^{-}\subset\Omega^{1}\Omega^{-}$
l'identité \[
\mathcal{D}^{2}=d^{*}d+dd^{*}-\frac{\scal}{12}.\]

\end{lem}
\begin{proof}
Sur $\Omega^{1}\Omega^{-}$, on a \begin{align*}
d^{*}d^{+}-d^{*}d^{-} & =-*d*(d^{+}-d^{-})=-*d^{2}=-*F(\Omega^{-})\\
d^{*}d^{+}+d^{*}d^{-} & =d^{*}d\end{align*}
 et par conséquent \[
d^{*}d^{-}=\frac{1}{2}(d^{*}d+*F(\Omega^{-})).\]
La courbure de $\Omega^{-}$ est réduite à la courbure scalaire :
\[
F(\Omega^{-})_{X,Y}=-\frac{\scal}{12}[X\wedge Y,\cdot].\]
Un calcul explicite facile donne \[
*F(\Omega^{-})|_{\mathcal{S}^{+}\mathcal{S}_{3}^{-}}=-\frac{\scal}{12}\]
d'où on déduit le lemme pour $\mathcal{D}^{2}=d^{*}d+2d^{*}d^{-}$. 
\end{proof}
\begin{cor}
\label{cor:imD}1) L'opérateur $\mathcal{D}^{2}$ sur $\mathcal{S}^{+}\mathcal{S}_{3}^{-}$
est un isomorphisme $H^{2}\rightarrow L^{2}$.

2) L'image par $\mathcal{D}$ de $H^{1}(\mathcal{S}^{+}\mathcal{S}_{3}^{-})$
dans $L^{2}(\mathcal{S}^{-}\mathcal{S}_{3}^{-})$ est fermée, égale
à $(\ker\mathcal{D})^{\perp}$.
\end{cor}
\begin{proof}
Le lemme \ref{lem:W} nous fournit une estimation \[
\|\mathcal{D}s\|^{2}\geq-\frac{\scal}{12}\| s\|^{2}\]
qui impose la première assertion (les opérateurs de ce type sont étudiés
dans \cite[I.2.B]{Biq00}).

La même estimation montre que l'image de $\mathcal{D}$ est fermée,
ce qui donne la seconde assertion.
\end{proof}

\subsubsection{Régularité dans les espaces à poids}

\begin{lem}
\label{lem:Dreg}L'opérateur $\mathcal{D}^{2}$ est un isomorphisme
\[
FS^{\ell+\delta}(\mathcal{U})\oplus H_{4+\delta}^{k+2;\ell}(\mathcal{S}^{+}\mathcal{S}_{3}^{-})\longrightarrow H_{4+\delta}^{k;\ell}(\mathcal{S}^{+}\mathcal{S}_{3}^{-}),\]
où $\mathcal{U}$ est le sous-fibré de $\mathcal{S}^{+}\mathcal{S}_{3}^{-}$
constitué des espaces propres pour les valeurs propres $\pm4i$ de
$\xi=\sigma_{1}^{+}+\sigma_{1}^{-}$. 
\end{lem}
\begin{proof}
Il s'agit d'une simple application de la théorie rappelée dans la
section \ref{sec:elliptique}. En effet, par le corollaire \ref{cor:imD},
l'opérateur $\mathcal{D}^{2}$ est inversible dans $L^{2}$. De plus,
sur $\mathcal{S}^{+}\mathcal{S}_{3}^{-}\subset\Omega^{1}\Omega^{-}$,
on a la formule de Weitzenböck \[
dd^{*}+d^{*}d=\nabla^{*}\nabla+\frac{\scal}{3},\]
qui, combinée au lemme \ref{lem:W}, fournit \[
\mathcal{D}^{2}=\nabla^{*}\nabla+\frac{\scal}{4}.\]
Le calcul des termes d'ordre 0 de $\nabla^{*}\nabla$ est simple :
\begin{align*}
-\sum_{1}^{3}\rho_{0}(Y_{i})^{2} & =-\sum_{1}^{3}\rho_{3}^{-}(\sigma_{i}^{-})^{2}-\left(\frac{1}{2}\rho_{3}^{-}(\sigma_{1}^{-})-\frac{3}{2}\rho^{+}(\sigma_{1}^{+})\right)^{2}+\rho_{3}^{-}(\sigma_{1}^{-})^{2}\\
 & =\fC(\rho_{3}^{-})+\frac{9}{4}+\frac{3}{4}\rho_{3}^{-}(\sigma_{1}^{-})^{2}+\frac{3}{2}\rho_{3}^{-}(\sigma_{1}^{-})\rho^{+}(\sigma_{1}^{+})\end{align*}
compte tenu de $\fC(\rho_{3}^{-})=15$, la plus petite valeur propre
est obtenue pour les valeurs propres $\pm(3i,i)$ de $(\rho_{3}^{-}(\sigma_{1}^{-}),\rho^{+}(\sigma_{1}^{+}))$,
elle est égale à $6$, donc la plus petite valeur propre des termes
d'ordre 0 de $\mathcal{D}^{2}$ est $\lambda=0$, ce qui fournit le
poids critique $4$ et le résultat du lemme.
\end{proof}

\subsubsection{Stabilité du noyau}

Revenons à l'opérateur de Dirac $\mathcal{D}$ sur $\mathcal{S}^{-}\mathcal{S}_{3}^{-}$
et à son noyau $L^{2}$. Perturbons la métrique de $\bC H^{2}$ en
une métrique ACH $g$ : puisque $\mathcal{D}^{2}$ est inversible
sur $\mathcal{S}^{+}\mathcal{S}_{3}^{-}$, l'opérateur $\mathcal{D}_{g}\mathcal{D}$
reste inversible pour une petite perturbation $g$, et on notera $G_{g}$
son inverse. Construisons alors l'opérateur $\Phi_{g}:\ker\mathcal{D}\rightarrow(\ker\mathcal{D})^{\perp}$
par la formule, pour $s\in L^{2}(\mathcal{S}^{-}\mathcal{S}_{3}^{-})$,
\[
\Phi_{g}(s)=-\mathcal{D}G_{g}\mathcal{D}_{g}s.\]
Alors $\mathcal{D}_{g}(s+\Phi_{g}(s))=0$, donc le noyau $L^{2}$
de $\mathcal{D}_{g}$ s'identifie au graphe de $\Phi_{g}$.

On en déduit la première assertion du lemme suivant.

\begin{lem}
\label{lem:Hfibre}L'application $s\rightarrow s+\Phi_{g}(s)$ identifie
$\ker\mathcal{D}$ et $\ker\mathcal{D}_{g}$ ; ainsi les noyaux des
opérateurs de Dirac sur $\mathcal{S}^{-}\mathcal{S}_{3}^{-}$ forment
un fibré hilbertien au-dessus de l'espace des métriques ACH (dans
un voisinage de $\bC H^{2}$).

Pour $g$ ayant la régularité décrite dans le lemme \ref{lem:regg},
cette application donne une identification \[
\Xi_{g}:(FS^{\ell+\delta}(\mathcal{J})\oplus H_{4+\delta}^{\infty;\ell})\cap\ker\mathcal{D}\longrightarrow(FS^{\ell+\delta}(\mathcal{J})\oplus H_{4+\delta}^{\infty;\ell})\cap\ker\mathcal{D}_{g},\]
 donc on obtient un fibré au-dessus de l'espace des métriques de cette
régularité.
\end{lem}
\begin{rem}
Dans le cas où la métrique $g$ est autoduale, le noyau de $\mathcal{D}_{g}$
est entièrement constitué de sections de $\mathcal{S}_{4}^{-}$. En
effet, d'après le lemme \ref{lem:DS4}, l'opérateur $\mathcal{D}_{g}^{2}$
préserve $\mathcal{S}_{4}^{-}$ et $\mathcal{S}_{2}^{-}$ ; mais sur
$\mathcal{S}_{2}^{-}=\Omega^{-}$, il s'identifie à $\nabla^{*}\nabla$
qui n'a pas de noyau.
\end{rem}
\begin{proof}
La seconde assertion mérite de vérifier la régularité des opérateurs.
Si la métrique a la régularité décrite, l'opérateur $\mathcal{D}_{g}\mathcal{D}$
est bien défini entre les espaces $FS^{\ell+\delta}(\mathcal{J})\oplus H_{4+\delta}^{k+1;\ell}\rightarrow H_{4+\delta}^{k-1;\ell}$.
Comme perturbation de l'opérateur $\mathcal{D}^{2}$ qui est un isomorphisme
par le lemme \ref{lem:Dreg}, il a un inverse $G_{g}$, donc l'application
$s\rightarrow s-\mathcal{D}G_{g}\mathcal{D}_{g}s$ fournit l'identification
voulue $(FS^{\ell+\delta}(\mathcal{J})\oplus H_{4+\delta}^{k;\ell})\cap\ker\mathcal{D}\rightarrow(FS^{\ell+\delta}(\mathcal{J})\oplus H_{4+\delta}^{k;\ell})\cap\ker\mathcal{D}_{g}$.
\end{proof}

\section{Métriques autoduales}

Nous sommes à présent en mesure de démontrer le résultat principal
de cet article.

\subsection{Schéma de la démonstration}

Notons $\Upsilon^{\ell+4+\delta}$ l'espace des structures CR de régularité
$FS^{\ell+4+\delta}$ sur $S^{3}$.

Une structure CR $J$ sur $S^{3}$, de régularité $FS^{\ell+4+\delta}$,
proche de $J_{0}$, est le bord d'une métrique ACH, Einstein, $g$,
de régularité précisée par le lemme \ref{lem:regg}. Pour toute métrique
d'Einstein, le tenseur de Weyl, en tant que 2-forme à valeurs dans
les 2-formes, est harmonique : \[
(d_{g}d_{g}^{*}+d_{g}^{*}d_{g})W_{g}^{-}=0.\]

\begin{lem}
Pour $g$ Einstein, dont l'infini conforme est de régularité $FS^{\ell+4+\delta}$,
le tenseur de Weyl $W_{g}^{-}$ satisfait \[
W_{g}^{-}=W_{4}^{-}e^{-4r}+O(e^{-(4+\delta)r}),\]
avec $W_{4}^{-}\in FS^{\ell+\delta}(\mathcal{J})$ et $W_{g}^{-}-W_{4}^{-}e^{-4r}\in H_{4+\delta}^{\infty;\ell}$.
\end{lem}
\begin{proof}
Par le lemme \ref{lem:regg}, la métrique s'écrit $g=\overline{g}+\Gamma e^{-2r}+G$,
avec $\Gamma\in FS^{\ell+\delta}$ et $G\in H_{4+\delta}^{\infty;\ell}$.
La métrique $\overline{g}$ est obtenue par un développement polyhomogène,
et par le lemme \ref{lem:Ricgbar2}, le tenseur $W_{\overline{g}}^{-}$
admet aussi un tel développement, avec premier terme non nul d'ordre
$4$, section du fibré $\mathcal{J}$ \cite[proposition 6.5]{BiqHer}
: il a donc la régularité annoncée dans le lemme.

Passons maintenant à la métrique $g$ : le tenseur $W_{g}^{-}$ a
encore un développement polyhomogène partiel (jusqu'à l'ordre $4$),
mais la perturbation apportée à $W_{\overline{g}}^{-}$ par le terme
$\Gamma e^{-2r}+G$ ne peut modifier que le terme d'ordre $4$, par
$\textrm{cst.}\Gamma e^{-2r}$ ; cela prouve le lemme.
\end{proof}
Par conséquent, l'application $g\rightarrow W_{g}^{-}$ constitue
une section du fibré des spineurs harmoniques, à bord de régularité
$FS^{\ell+\delta}$. Il sera commode, via l'identification $\Xi_{g}$
construite dans le lemme \ref{lem:Hfibre}, de considérer plutôt \[
g\longrightarrow\Xi_{g}^{-1}W_{g}^{-}\]
qui est à valeurs dans l'espace vectoriel fixe $(FS^{\ell+\delta}(\mathcal{J})\oplus H_{4+\delta}^{\infty;\ell})\cap\ker\mathcal{D}_{g_{0}}$.

On peut regarder l'application entre les valeurs à l'infini : à une
structure CR $J$ on associe la valeur à l'infini $\partial_{\infty}W_{g}^{-}$
du spineur harmonique $W_{g}^{-}$, où $g$ est la métrique d'Einstein
qui remplit $J$. Ramenant comme ci-dessus les spineurs à la métrique
fixe $g_{0}$, on obtient un opérateur $W_{\infty}^{-}(J)=\partial_{\infty}\Xi_{g}^{-1}W_{g}^{-}$,
défini entre les espaces\[
W_{\infty}^{-}:\Upsilon^{\ell+4+\delta}\longrightarrow FS^{\ell+\delta}(\mathcal{J}^{+}).\]
Notons que le choix de jauge sur la métrique d'Einstein $g$ qui remplit
$J$ n'influe pas sur $W_{\infty}^{-}$, car les difféomorphismes
égaux à l'identité sur le bord à l'infini ne modifient pas $W_{\infty}^{-}(J)$.

Le lemme principal duquel nous déduirons le théorème est le suivant
:

\begin{lem}
\label{lem:surj}L'opérateur $d_{g_{0}}W_{\infty}^{-}$, restreint
aux déformations remplissables par une métrique Kähler-Einstein, est
surjectif.
\end{lem}
Autrement dit : tous les spineurs harmoniques dans $\mathcal{S}_{4}^{-}$
sont obtenus comme des tenseurs de Weyl de métriques Kähler-Einstein
infinitésimales. Ce lemme sera démontré dans la section \ref{sec:33}. 

\begin{rem}
Comme conséquence du lemme, on a une manière (compliquée) de montrer
que toutes les sections de $\mathcal{J}^{+}$ de régularité $FS^{\ell+\delta}$
sont des bords de spineurs harmoniques dont on connaît la régularité,
voir remarque \ref{rem:continuite}.
\end{rem}
On en déduit le résultat suivant, qui implique immédiatement le théorème
principal.

\begin{lem}
\label{lem:transverse}1) L'opérateur $W_{\infty}^{-}$ est submersif.

2) Le noyau de la différentielle $d_{g_{0}}W_{\infty}^{-}$ en la
métrique hyperbolique complexe est transverse aux structures CR qui
sont les bords de métriques Kähler-Einstein, et l'intersection est
réduite à l'image infinitésimale des contactomorphismes.
\end{lem}
\begin{proof}
La première assertion résulte directement du lemme \ref{lem:surj}.
On en déduit que l'espace des $J\in\Upsilon^{\ell+\delta}$ qui sont
le bord d'une métrique ACH autoduale d'Einstein forment une sous-variété
de $\Upsilon^{\ell+\delta}$, d'espace tangent en $g_{0}$ le noyau
de $d_{g_{0}}W_{\infty}^{-}$.

L'assertion sur la transversalité est aussi une conséquence directe
du lemme \ref{lem:surj}. Il en résulte que l'intersection représente
l'espace tangent aux structures CR $J$ dont le remplissage Einstein
est à la fois autodual et Kähler : cette intersection peut être calculée
explicitement, voir la section \ref{sec:33}. Cependant il y a un
raisonnement géométrique plus direct : on sait que la seule métrique
autoduale et Kähler-Einstein est la métrique hyperbolique complexe,
donc $J$ ne peut être que la structure CR standard, à un contactomorphisme
près.
\end{proof}

\subsection{Bords de métriques Kähler-Einstein}

Un théorème fondamental, dû à Burns et Epstein, Lempert, Bland \cite{BurEps90,Eps92,Lem92,Bla94}
assure que la structure CR $J$ sur $S^{3}$ est remplissable par
une structure complexe si et seulement si, à l'action près d'un contactomorphisme,
le tenseur $\phi$ correspondant via (\ref{eq:phi}) n'a que des coefficients
de Fourier positifs ou nuls par rapport à l'action de $U_{1}$ sur
$S^{3}$.

De plus, Bland a construit une forme normale pour une structure CR
$J$, par rapport à l'action des contactomorphismes, dans laquelle
$J$ est remplissable si et seulement si ses coefficients de Fourier
négatifs s'annulent. Infinitésimalement, cela donne le théorème suivant
(rappelons la notation $S_{4}^{-}(\pm4)$ provenant de la définition
\ref{def:Jplus}) :

\begin{thm}
[Bland]Une structure CR infinitésimale $\dot{J}\in FS^{\ell}(\mathcal{J})$
est le bord d'une déformation complexe de la boule si et seulement
si ses composantes dans $V_{\rho}\otimes(V_{\rho}S_{4}^{-}(\pm4))^{S^{1}}$,
où $V_{\rho}=\bC_{K}S_{L}^{-}$, sont nulles dès que $|K|>L$.
\end{thm}
\begin{notation}
On notera $KE^{\ell}(\mathcal{J})$ l'espace des déformations infinitésimales
de $J_{0}$ qui sont remplissables par une métrique Kähler-Einstein.
Par le théorème de Bland, cet espace est donné par \[
KE^{\ell}(\mathcal{J})=FS^{\ell}(\mathcal{J})\cap\oplus_{|K|\leq L}V_{\rho}\otimes(V_{\rho}S_{4}^{-}(\pm4))^{S^{1}}.\]

\end{notation}
Pour calculer la différentielle du tenseur de Weyl, nous avons besoin
de donner un peu plus de détail sur le théorème de Bland.

Rappelons qu'une structure CR sur $S^{3}$ est fournie par une section
$\phi$ de $\Omega^{0,1}\otimes T^{1,0}$, sur lequel $\xi=\sigma_{1}^{+}+\sigma_{1}^{-}$
agit avec poids $4$. La décomposition harmonique est donc \[
L^{2}(\Omega^{0,1}\otimes T^{1,0})=\oplus_{\rho}V_{\rho}\otimes(V_{\rho}\otimes\bC_{4})^{S^{1}}.\]
En particulier, pour $V_{\rho}=\bC_{K}S_{L}^{-}$, un $w\in(V_{\rho}\otimes\bC_{4})^{S^{1}}$
admet un poids \begin{equation}
k=-K-4\label{eq:kphi}\end{equation}
pour l'action de $\sigma_{1}^{-}$. En particulier, $V_{\rho}\otimes(V_{\rho}\otimes\bC_{4})^{S^{1}}$
est non nul seulement pour \begin{equation}
-L-4\leq K\leq L-4.\label{eq:Kphicon}\end{equation}

D'autre part, toujours d'après Bland, les contactomorphismes infinitésimaux,
de régularité $FS^{\ell+1}$, sont paramétrés par une fonction réelle
$f\in FS^{\ell+2}$, de sorte que l'action infinitésimale des contactomorphismes
soit \begin{equation}
f\longrightarrow\overline{\partial}\sharp\overline{\partial}f,\label{eq:act}\end{equation}
où $\sharp:\Omega^{0,1}\rightarrow T^{1,0}$ est l'identification
induite par la forme symplectique $d\eta$. Au niveau des décompositions
harmoniques, la fonction $f$ est représentée par des termes \[
v\otimes w\in V_{\rho}\otimes V_{\rho}^{S^{1}},\]
et le poids de $\sigma_{1}^{-}$ est \begin{equation}
k=-K,\label{eq:kf}\end{equation}
de sorte que $K$ est soumis à la contrainte \begin{equation}
-L\leq K\leq L.\label{eq:Kfcon}\end{equation}
La condition de réalité se traduit par l'invariance sous \[
v\otimes w\longrightarrow(\tau v)\otimes(\tau\overline{w}),\]
où $\tau$ est une structure réelle sur $S_{L}^{-}$, donc \[
\tau:\bC_{K}S_{L}^{-}\longrightarrow\bC_{-K}S_{L}^{-}.\]
L'entier $K$ représente le poids de l'action de $U_{1}$ : cette
condition dit seulement que les coefficients à $K>0$ de la fonction
réelle sont déterminés par ceux à $K<0$.

L'action infinitésimale des contactomorphismes (\ref{eq:act}) s'écrit
\[
w\longrightarrow-i\rho(Y)^{2}w,\]
où $(H,X,Y)$ est la représentation standard de $\mathfrak{sl}_{2}$
associée à $(\sigma_{i})$, de sorte que $[H,X]=2X$, $[H,Y]=-2Y$
et $[X,Y]=H$. Cette action diminue bien le poids $k$ de $4$, comme
exigé par les égalités (\ref{eq:kphi}) et (\ref{eq:kf}).

La jauge de Bland consiste à faire agir les contactomorphismes de
sorte de tuer, autant que possible, les coefficients à $K<0$ dans
$\phi$. En vue de (\ref{eq:Kphicon}) et (\ref{eq:Kfcon}), on arrive
à tuer tous ces coefficients, sauf ceux d'ordre $K=-L-4$ ou $-L-2$.
Compte tenu de la correspondance $\dot{J}=\phi+\overline{\phi}$,
cela donne le théorème énoncé plus haut.

\subsection{Calcul du tenseur de Weyl\label{sec:33}}

\subsubsection{Cas des métriques Kähler-Einstein}

Dans le cas où $J$ est l'infini conforme d'une métrique Kähler-Einstein
$g$, on peut calculer $\partial_{\infty}W_{g}^{-}$ : 

\begin{lem}
[\protect{\cite[proposition 6.5]{BiqHer}}]\label{lem:WQ}Pour $g$
Kähler-Einstein, d'infini conforme $J$, la valeur à l'infini du tenseur
de Weyl est \[
\partial_{\infty}W_{g}^{-}=\textrm{cst. }Q(J),\]
où $Q$ est la courbure de Cartan de la structure $J$.
\end{lem}
La différentielle en $J_{0}$ de $W_{\infty}^{-}$ fournit un opérateur
\[
d_{g_{0}}W_{\infty}^{-}:FS^{\ell+4+\delta}(\mathcal{J})\longrightarrow FS^{\ell+\delta}(\mathcal{J}^{+})\]
qui à chaque déformation infinitésimale $\dot{J}$ associe la valeur
à l'infini du spineur harmonique fourni par le $W^{-}$ du remplissage
d'Einstein. Si on se restreint aux déformations Kähler-Einstein, le
lemme précédent indique qu'on a seulement la différentielle du tenseur
de Cartan, qui a été étudiée dans \cite{CheLee90} :

\begin{thm}
[Cheng et Lee]\label{thm:CL}L'opérateur $d_{J_{0}}Q$ est un opérateur
hypoelliptique d'ordre $4$, transversalement à l'action des contactomorphismes
de $S^{3}$, et son noyau est réduit à l'action infinitésimale des
contactomorphismes.
\end{thm}
\begin{rem}
Cheng et Lee montrent aussi que l'image de $d_{J_{0}}Q$ est le noyau
d'un opérateur de Bianchi, d'ordre 2 (ce qui est la contrepartie de
l'invariance de $d_{J_{0}}Q$ sous les contactomorphismes). Du point
de vue des spineurs harmoniques sur $\bC H^{2}$, il est plausible
que cet opérateur de Bianchi soit exactement la contrainte sur les
valeurs à l'infini, mentionnée à la remarque \ref{rem:obstruction}.
\end{rem}

\subsubsection{Démonstration du lemme \ref{lem:surj}}

Vu le lemme \ref{lem:WQ}, il s'agit de montrer que \[
d_{J_{0}}Q:KE^{\ell+4+\delta}(\mathcal{J})\longrightarrow FS^{\ell+\delta}(\mathcal{J}^{+})\subset FS^{\ell+\delta}(\mathcal{J})\]
 est surjectif.

D'après le théorème \ref{thm:CL}, l'image de $d_{J_{0}}Q$ est fermée
(en fait égale au noyau d'un opérateur de Bianchi), donc il suffit
de tester que $d_{J_{0}}Q$ est surjectif pour chaque représentation
$\rho$ concernée. Par le même théorème, le noyau de $d_{J_{0}}Q$
est égal à l'image infinitésimale des contactomorphismes, donc la
vérification de la surjectivité est réduite à un simple compte de
dimensions.

Continuons donc le calcul dans la décomposition harmonique, commencé
dans la section précédente. Pour chaque représentation $\rho$, l'opérateur
est à valeurs dans $J_{\rho}^{+}$ (définition \ref{def:Jplus}),
qui est non nul seulement pour $|K|\leq L-4$, et de dimension égale
à $\dim V_{\rho}$.

Le seul cas à considérer est donc $|K|\leq L-4$ : comme les objets
à considérer sont réels, nous regardons simultanément les représentations
$\bC_{K}S_{L}^{-}$ et $\bC_{-K}S_{L}^{-}$. Les structures CR infinitésimales
y forment un espace de dimension réelle $4\dim S_{L}$, les contactomorphismes
infinitésimaux un espace de dimension réelle $2\dim S_{L}$, et les
spineurs harmoniques un espace de dimension réelle $2\dim V_{\rho}$.
L'opérateur $d_{J_{0}}Q$ est donc surjectif.\qed

\begin{rem}
Pour $|K|=L-2$ ou $L$, l'opérateur $d_{J_{0}}Q$ s'annule, et on
peut vérifier effectivement que les déformations correspondantes sont
dans l'image infinitésimale des contactomorphismes.
\end{rem}

\subsubsection{Calcul de l'espace tangent aux bords de métriques autoduales}

Finalement, on peut déduire du calcul la précision suivante sur les
structures CR qui sont des infinis conformes de métriques autoduales
Einstein :

\begin{thm}
\label{thm:tangent}L'espace tangent en $J_{0}$ aux structures CR
de $\Upsilon^{\ell+4+\delta}$ qui sont remplissables par une métrique
autoduale d'Einstein est égal, modulo l'action infinitésimale des
contactomorphismes, à \[
FS^{\ell+4+\delta}(\mathcal{J})\cap\bigoplus_{|K|=L+2,L+4}V_{\rho}\otimes(V_{\rho}S_{4}^{-}(\pm4))^{S^{1}}.\]
\qed
\end{thm}
\bibliographystyle{smfalpha}
\bibliography{biblio,biquard}

\end{document}